\documentclass[MSNbibl,number,citesort,dvips]{arxbj}
\usepackage{graphicx}
\aid{0}
\volume{19}
\issue{5A}
\pubyear{2013}
\firstpage{2067}
\lastpage{2097}
\doi{10.3150/12-BEJ443} %kopijuoti is PTS

\makeatletter
\def\eqref#1{(\ref{#1})}

\newcommand{\rrVert}{\Vert}
\newcommand{\rrvert}{\vert}
\newcommand{\llVert}{\Vert}
\newcommand{\llvert}{\vert}
\newtheorem{theo}{Theorem}
\newtheorem{coro}{Corollary}
\makeatother

\begin{document}
\begin{frontmatter}

\title{Confidence bands for Horvitz--Thompson estimators using sampled
noisy functional~data}
\runtitle{Confidence bands for Horvitz--Thompson estimators}

\begin{aug}
%%%% inicialai - be tarpu
\author[1]{\fnms{Herv\'e} \snm{Cardot}\corref{}\thanksref{1,e1}\ead[label=e1,mark]{herve.cardot@u-bourgogne.fr}},
\author[2]{\fnms{David} \snm{Degras}\thanksref{2}\ead[label=e2]{ddegrasv@depaul.edu}} \and
\author[1]{\fnms{Etienne} \snm{Josserand}\thanksref{1,e3}\ead[label=e3,mark]{etienne.josserand@u-bourgogne.fr}}
\runauthor{H. Cardot, D. Degras and E. Josserand} %% auto
\address[1]{Institut de Math\'ematiques de Bourgogne, UMR 5584,
Universit\'e de Bourgogne, 9 Avenue Alain Savary, 21078 Dijon,
France.\\
\printead{e1,e3}}
\address[2]{DePaul University, 2320 N. Kenmore Avenue, Chicago, IL
60614, USA.\\ \printead{e2}}
\end{aug}

% HISTORY:
\received{\smonth{5} \syear{2011}}
\revised{\smonth{3} \syear{2012}}

% ABSTRACT
%
\begin{abstract}
When collections of functional data are too large to
be exhaustively observed,
survey sampling techniques
provide an effective way to estimate global quantities such as the
population mean function.
Assuming functional data are collected from a finite population
according to a probabilistic sampling scheme,
with the measurements being discrete in time and noisy,
we propose to first smooth the sampled
trajectories with local polynomials and then estimate the mean function
with a Horvitz--Thompson estimator.
Under mild conditions on the population size, observation times,
regularity of the trajectories, sampling scheme, and smoothing bandwidth,
we prove a Central Limit theorem in the space of continuous functions.
We also establish the uniform consistency of a covariance function estimator
and apply the former results to build confidence bands for the mean function.
The bands attain nominal coverage and are obtained through Gaussian
process simulations conditional on the estimated covariance function.
To select the bandwidth, we propose a cross-validation method that
accounts for the sampling weights.
A simulation study assesses the performance of our approach and
highlights the influence of the sampling scheme and bandwidth choice.
\end{abstract}

% KEYWORDS
%

\end{frontmatter}
%

%s1 ###
\section{Introduction}\label{sec1}

The recent development of automated sensors has given access to very
large collections of signals sampled at fine time scales.
However, exhaustive transmission, storage, and analysis of such massive
functional data may incur very large investments.
In this context, when the goal is to assess a global indicator like the
mean temporal signal, survey sampling techniques
are appealing solutions as they offer a good trade-off between
statistical accuracy and global cost of the analysis.
In particular, they are competitive with signal compression techniques
(Chiky and H\'ebrail \cite{ChHe08}).
The previous facts provide some explanation why, although survey
sampling and functional data analysis have been long-established
statistical fields,
motivation for studying them jointly only recently emerged in the literature.
In this regard, Cardot \textit{et al.} \cite{CaChGoLa10a} examine the theoretical
properties of functional principal components analysis (FPCA) in the
survey sampling framework.
Cardot \textit{et al.} \cite{CaDeJo10b} harness FPCA for model-assisted estimation
by relating the unobserved principal component scores to available
auxiliary information.
Focusing on sampling schemes, Cardot and Josserand \cite{CaJo11} estimate the
mean electricity consumption curve in a population of about 19,000
customers whose electricity meters were read every 30 minutes during
one week.
Assuming exact measurements, they first perform a linear interpolation
of the discretized signals and then consider a functional version of
the Horvitz--Thompson estimator.
For a fixed sample size, they show that estimation can be greatly
improved by utilizing stratified sampling over simple random sampling
and they extend the Neyman optimal allocation rule (see, e.g., S\"
{a}rndal \textit{et al.} \cite{SaSwWr92}) to the functional setup.
Note however that the finite-sample and asymptotic properties of their
estimator rely heavily on the assumption of error-free measurements,
which is not always realistic in practice.
%
%it might be good to pull the contribution up, follow it with the lit
%review and then give a brief summary. as it stands now, i'm missing
%the reason why i'm told all the details about herve's recent papers.
%if you start with "the first contribution ... is to generalize the
%framework of c&j and then give the literature review, it would be much
%clearer.
The first contribution of the present work is to generalize the
framework of Cardot and Josserand \cite{CaJo11} to noisy functional data.
Assuming curve data are observed with errors that may be correlated
over time,
we replace the interpolation step in their procedure by a smoothing
step based on local polynomials.
As opposed to interpolation, smoothing can effectively reduce the noise
level in the data, which improves estimation accuracy.
%We extend the previous asymptotic theory by establishing a functional
%CLT for the resulting mean function estimator
%and by proving the uniform consistency of a related covariance
%estimator.
We establish a functional CLT for the mean function estimator based on
the smoothed data and prove the uniform consistency of a related
covariance estimator. These results have important applications to the
simultaneous inference of the mean function.
%can you give a brief preview what great thing this accomplishes?

In relation to mean function estimation, a key statistical task is to
build confidence regions.
%could you make the relation to *your* project a little clearer? again,
%we are going into a very extensive literature review and i'm not sure
%where we are going. why do i need to know all that?
There exists a vast and still active literature on confidence bands in
nonparametric regression.
See, for example, Sun and Loader \cite{SuLo94}, Eubank and Speckman \cite{EuSp93},
Claeskens and van Keilegom \cite{ClKe03}, Krivobokova \textit{et al.} \cite{KrKnCl10},
and the references therein.
%However,
When data are functional the~literature is much less abundant.
%can you give us a hint as to *why* that is? this could help you
%establish how important your work is.
One possible approach is to obtain confidence balls for the mean
function in a $L^2$-space.
Mas \cite{Ma07} exploits this idea in a goodness-of-fit test based on the
functional sample mean and regularized inverse covariance operator.
Using adaptive projection estimators, Bunea \textit{et al.} \cite{BuIvWe11}
build conservative confidence regions for the mean of a Gaussian process.
Another approach consists in deriving results in a space $C$ of
continuous functions equipped with the supremum norm.
This allows for the construction of confidence bands that can easily be
visualized and interpreted, as opposed to $L^2$-confidence balls.
This approach is adopted, for example, by Faraway \cite{Fa97} to build
bootstrap bands in a varying-coefficients model, by Cuevas \textit{et al.}
\cite{CuFeFr06} to derive bootstrap bands for functional location parameters, by
Degras \cite{De09,De11} to obtain normal and bootstrap bands using noisy
functional data, and by Cardot and Josserand \cite{CaJo11} in the context of a
finite population.
In the latter work, the strategy was to first establish a CLT in the
space $C$ and then derive confidence bands based on a simple but rough
approximation to the supremum of a Gaussian process (Landau and Shepp \cite{LaSh70}).
Unfortunately, the associated bands depend on the data-generating
process only through its variance structure and not its correlation
structure, which may cause the empirical coverage to differ from the
nominal level.
% again, pull up the main claim! point first...
The second innovation of our paper is to propose confidence bands that
are easy to implement
and attain nominal coverage in the survey sampling/finite population setting.
To do so, we use Gaussian process simulations as
in Cuevas \textit{et al.} \cite{CuFeFr06} or Degras \cite{De11}. This procedure can
be thought as a parametric bootstrap, where the parameter to be
estimated, the covariance function, is lying in an infinite dimensional
functional space.
Our contribution is to provide the theoretical underpinning of the
construction method,
thereby guaranteeing that nominal coverage is attained asymptotically.
The theory we derive involves maximal inequalities,
random entropy numbers, and large covariance matrix theory.

%, which can be viewed as a parametric bootstrap in an infinite
%dimension context, is much more effective and consists in estimating
%the covariance function of the in order to
% of this process and determining the quantiles of the supremum.

Finally, the implementation of the mean function estimator developed in
this paper requires the selection of a bandwidth in the data smoothing step.
Objective, data-driven bandwidth selection methods are desirable for
this purpose.
%again, i would love to know why this is important.
As explained by Opsomer and Miller \cite{OpMi05}, bandwidth selection in the
survey estimation context poses specific problems (in particular, the
necessity to take the sampling design into account) that make usual
cross-validation or mean square error optimization methods inadequate.
In view of the model-assisted survey estimation of a population total,
these authors propose a cross-validation method that aims at minimizing
the variance of the estimator, the bias component being negligible in
their setting. In our functional and design-based framework, the bias
is however no longer negligible.
%pull this up! it would be soooo much stronger to hear about your
%contribution first
We therefore devise a novel cross-validation criterion based on
weighted least squares, with weights proportional to the sampling weights.
For the particular case of simple random sampling without replacement,
this criterion reduces to the cross validation technique of Rice and
Silverman \cite{RiSi91}, whose asymptotic properties has been studied by Hart
and Wehrly \cite{HaWe93}.
%"... and our paper is therefore incredibly important for that-and-that
%reason."

%For unequal probability sampling, Berger and Skinner (2005) and % have
%already been proposed for by. or jackknife

The paper is organized as follows. We fix notations and define our
estimators in Section~\ref{sec2}.
In Section~\ref{sec3}, we introduce our asymptotic framework based on
superpopulation models (see Isaki and Fuller \cite{IsFu82}),
establish a CLT for the mean function estimator %of the mean trajectory
in the space of continuous functions, and show the uniform consistency
of a covariance estimator.
%After that, we prove that by simulating the limiting Gaussian process
%conditional on its estimated covariance,
%one can build confidence bands that have asymptotically correct
%coverage.
Based on these results, we propose a simple and effective method for
building simultaneous confidence bands.
%The bands only rely on Gaussian process simulations and have
%asymptotic nominal coverage. }
In Section \ref{sec4}, a weighted cross-validation procedure is proposed for
selecting the bandwidth
and simulations are performed to compare different sampling schemes and
bandwidth choices.
Our estimation methodology is seen to compare favorably with other
methods and to achieve nearly optimal performances.
% to assess the numerical performance of.
%and are useful to highlight the role of the bandwidth and the
%performances of the suggested cross-validation criterion.
The paper ends with a short discussion on topics for future research.
Proofs are gathered in an \hyperref[app]{Appendix}.

%%%%%%%%%%%%%%%%%%%%%%%%%%%%%%%%%%%%%%%%%%%%%%%
%%%%%% Introduction %%%%%
%%%%%%%%%%%%%%%%%%%%%%%%%%%%%%%%%%%%%%%%%%%%%%%
%s2 ###
\section{Notations and estimators}\label{sec2}

Consider a finite population $U_N = \{1,\dots,N\}$ of size $N$
and suppose that to each unit $k \in U_N$ corresponds a real function
$X_k$ on $ [0,T],$ with $T<\infty.$
We assume that each trajectory $X_k$ belongs to the space of continuous
functions $C([0,T]).$
Our target is the mean trajectory $\mu_N(t),  t \in[0,T],$ defined
as follows:
\begin{eqnarray}\label{def:mu}
\mu_N(t) &=& \frac{1}{N} \sum_{k \in U}
X_k(t) . % ,   t \in
%[0,T].
\end{eqnarray}

We consider a random sample $s$ drawn from $U_N$ without replacement
according to a fixed-size sampling design $p_N(s),$ where $p_N(s)$ is
the probability of drawing the sample $s.$ The size $n_N$ of $s$ is
nonrandom %and is equal to
and we suppose that the first and second order inclusion probabilities satisfy
%%%%% General case
%
\begin{itemize}
\item$\pi_k := \mathbb{P}(k \in s) >0$ for all $k \in U_N$
\item$\pi_{kl} := \mathbb{P}(k \& l \in s)>0$ for all $k,l \in U_N$
%%, $k
\end{itemize}
so that each unit and each pair of units can be drawn with a non null
probability from the population.
Note that for simplicity of notation the subscript $N$ has been omitted.
Also, by convention, we write $\pi_{kk} = \pi_k$ for all $k\in U_N$.

Assume that noisy measurements of the sampled curves are available
at $d=d_N$ fixed discretization points $0=t_1<t_2 < \cdots<t_d=T.$ For
all units $k \in s$, we observe %the measured data are written as
\begin{equation}
Y_{jk} = X_k(t_j) + \varepsilon_{jk},
\end{equation}
%
%where $\varepsilon_{jk}$ are realizations (not necessarily independent)
%of a centered random variable with finite variance, $\E\varepsilon_{jk}
%=0$ and $\E\varepsilon_{jk}^2 < C.$ We suppose that the selection
%process $p(s)$ of the sample $s$ is independent of the noise $
%random.
where the measurement errors $\varepsilon_{jk}$ are centered random
variables that are independent across the index $k$ (units) but not
necessarily across $j$ (possible temporal dependence).
%In other words, the as being independent across units but possibly
%correlated over time.
%We assume that the $\varepsilon_{jk}$ have uniformly bounded covariance
%structures, i.e. denoting by $\mathbf{V}_k$ the $p\times p$ covariance
%matrix of the random vector $(\varepsilon_{k1},\ldots,\varepsilon_{kp})$, it
%holds that $\| \mathbf{V}_k \| < C$ for all $k\in U_N$ and $N$, where $
It is also assumed that the random %selection process $p_N(s)$ of the
sample $s$ is independent of the noise $\varepsilon_{jk}$ and the
trajectories $X_k(t), t \in[0,T]$ are deterministic.

Our goal is to estimate $\mu_N$ as accurately as possible and to build
asymptotic confidence bands, as in Degras \cite{De11} and Cardot and
Josserand \cite{CaJo11}. For this, we must have a uniformly consistent
estimator of its covariance function.

%s2.1 ###
\subsection{Linear smoothers and the Horvitz--Thompson estimator}\label{sec2.1}

For each (potentially observed)
unit $k \in U_N$, we aim at recovering the curve $X_k$ by smoothing the
corresponding discretized trajectory $(Y_{1k}, \ldots, Y_{dk})$
with a linear smoother (e.g., spline, kernel, or local polynomial):
\begin{equation}\label{def:smoother}
\widehat{X}_k(t) = \sum_{j=1}^d
W_j(t)Y_{jk}. %  t \in[0,T].
\end{equation}
Note that the reconstruction can only be performed for the observed
units $k\in s$.

% and are supposed to satisfy (22a) and (22b) in Lemma 1 of Degras
%(2009, preprint).

Here we use local linear smoothers (see, e.g., Fan and
Gijbels \cite{f97})
because of their wide popularity, good statistical properties, and
mathematical convenience.
The weight functions $W_j(t)$ can be expressed as
\begin{equation}
\label{local linear 1} W_j(t) = \frac{ ({1}/{(dh)}) \{ s_{2}(t) - (t_j-t) s_{1}(t)
\} K  ( {(t_j-t)}/{h}  )}{ s_2(t)s_0(t) -s_1^2(t)}   ,\qquad    j=1,
\ldots,d,
\end{equation}
where $K$ is a kernel function, $h>0$ is a bandwidth, and
\begin{equation}
s_{l}(x)=  \frac{1}{dh} \sum
_{j=1}^d (t_j-t)^l   K
\biggl( \frac{t_j-t }{h} \biggr), \qquad   l=0,1,2.
\end{equation}
We suppose that the kernel $K$ is nonnegative, has compact support,
satisfies $K(0)>0$ and $|K(s)-K(t)| \le C |s-t|$ for some finite
constant $C$ and
for all $s,t \in[0,T]$.

The classical Horvitz--Thompson estimator of the mean curve is
\begin{eqnarray}\label{def:smoothmean}
\widehat{\mu}_N(t) &=& \frac{1}{N} \sum
_{k \in s} \frac{\widehat
{X}_k(t)}{\pi_k}
\nonumber\\[-8pt]\\[-8pt]
& = & \frac{1}{N} \sum_{k \in U}
\frac{\widehat{X}_k(t)}{\pi_k} I_k, %  t \in[0,T],
\nonumber
\end{eqnarray}
where $I_k$ is the sample membership indicator ($I_k = 1$ if $k \in s$
and $I_k = 0$ otherwise).
It holds that $\mathbb{E}(I_k) = \pi_k$ and $\mathbb{E}(I_k I_l) =
\pi_{kl}$.

%s2.2 ###
\subsection{Covariance estimation}\label{sec2.2}

The covariance function of $\widehat{\mu}_N$ can be written as
\begin{equation}
\label{cov expression} \operatorname{Cov} \bigl(\widehat{\mu}_N(s),
\widehat{\mu}_N(t) \bigr) = \frac{1}{N}\gamma_N(s,t)
\end{equation}
for all $s,t\in[0,T]$, where
\begin{equation}
\label{cov expression 2} %&= &\E\left( \left(\widehat{\mu}_N(t) - \E( \widehat{\mu}_N(t) )
\gamma_N(s,t) = \frac{1}{N} \sum_{k,l \in U}
\Delta_{kl} \frac
{\tilde{X}_k(s)}{\pi_k} \frac{\tilde{X}_l(t)}{\pi_l} + \frac
{1}{N}
\sum_{k \in U} \frac{1}{\pi_k} \mathbb{E} \bigl( \tilde
{\varepsilon }_k(s) \tilde{\varepsilon}_k(t) \bigr)
\end{equation}
with
\begin{equation}
\label{def tilde X eps Delta} \cases{\displaystyle
\tilde{X}_k(t) = \sum_{j=1}^d
W_j(t)X_k(t_j),
\cr
\displaystyle\tilde{\varepsilon}_k(t) = \sum_{j=1}^d
W_j(t)\varepsilon_{kj},
\cr
\Delta_{kl} = \operatorname{Cov}(I_k,I_l) =
\pi_{kl} - \pi_k\pi_l . }
\end{equation}

%, i.e. $\Delta_{kl}=\pi_{kl} - \pi_k \pi_l$ if $k \neq l$ and $

A natural estimator of $\gamma_N(s,t)$ is given by
%the covariance of $\widehat{\mu}_N(s,t)$ (up to the factor $1/N$),
%we use the unbiased estimator
%
\begin{equation}
\label{covariance estimator} \widehat{\gamma}_N(s,t) =
\frac{1}{N} \sum_{k,l\in U} \frac{\Delta_{kl}}{\pi_{kl}}
\biggl( \frac{I_k}{\pi_k}   \frac{I_l}{\pi_l} \biggr) \widehat{X}_{k}(s)
\widehat{X}_{l}(t).
\end{equation}
It is unbiased and its uniform mean square consistency is established
in Section \ref{consistency in covariance estimation}.
%Note that this estimator is unbiased and that it is a covariance
%function in the case of simple random sampling and stratified sampling.

%%%%%%%%%%%%%%%%%%%%%%%%%%%%%%%%%%%%%%%%%%%%%%%%%%%%%%%%%%%%%%%

%s3 ###
\section{Asymptotic theory}\label{sec3}

We consider the superpopulation framework introduced by Isaki and
Fuller \cite{IsFu82} and discussed in detail by Fuller \cite{Fu09}.
Specifically, we study the behaviour of the estimators $ \widehat{\mu
}_N$ and $ \widehat{\gamma}_N$ as
population $U_N = \{ 1,\ldots, N\}$ increases to infinity with $N$.
Recall that the sample size $n$, inclusion probabilities $\pi_k$ and
$\pi_{kl}$, and grid size $d$ all depend on $N$.
In what follows, we use the notations $c$ and $C$ for finite, positive
constants whose value may vary from place to place.
%(Formally we consider the limit in the increasing sequence of sets
%$U_1 \subset\cdots\subset U_N \subset\cdots$)
%size $N$ tends to infinity.
%Recall that for all $N\ge1$, a finite population of size $N$ is
%represented by the set $U_N= \{ 1,\ldots,N\}$.
%Now, for all $N\ge1$, let $s_N$ be a random sample of size $n_N$
%drawn from $U_N$ without replacement, according to the fixed-size
%sampling design $p_N (s_N )$.
%For simplicity of notation, we drop the subscript $N$ in the elements
%$s_N, p_N, \pi_{kN}, \pi_{klN}$ defined in Section \ref{sec: notation
%and estimators}.
% by $\pi_{kN}$ and $\pi_{klN}$ their first and second order inclusion
%probabilities.
%Note that the sequence of sub populations is an increasing nested one
%while the sample sequence is not
%in the following when there is no ambiguity.
The following assumptions are needed for our asymptotic
study.\looseness=-1
\begin{enumerate}[(A4)]

\item[(A1)] (\textit{Sampling design})
$ \frac{n}{N} \ge c ,   \pi_k \ge c ,  \pi_{kl} \ge
c , $ and
$   n | \pi_{kl} - \pi_k \pi_l | \le C $ for all $k, l
\in U_N$ ($k\ne l$) and $ N\ge1$.

%$ \liminf_{N\to\infty} (n/N) >0 $,
%$  \liminf_{N\to\infty} \min_{k\in U_N} \pi_k >0, $ $
% \liminf_{N\to\infty} \min_{ k \ne l} \pi_{kl} >0, $\\
%and
% $  \limsup_{N\to\infty}  n \max_{k\ne l} | \pi_{kl} -

\item[(A2)] (\textit{Trajectories})
$| X_k(s) - X_k(t) | \le C |s-t|^{\beta}$ and $ |X_k(0)| \le C $
for all $k \in U_N,  N\ge1,$ and $s,t\in[0,T]$, where $\beta> \frac
{1}{2}$ is a finite constant.

\item[(A3)] (\textit{Growth rates}) % of sample size and discretization
%grid})
$c \le d (t_{j+1}-t_j) \le C$ for all $1\le j \le d,   N\ge1,$
and %$  \lim_{N\to\infty} (N/d) = C$, where $C>0$ is a
%finite constant.
$  \frac{d(\log\log N) }{N} \to0$  as $N\to\infty$.
%$  \lim_{N\to\infty} N \left( \frac{\log d}{d}\right)^{2

\item[(A4)] (\textit{Measurement errors}) The random vectors $(\varepsilon_{k1},\ldots,\varepsilon_{kd})', k\in U_N,$ are i.i.d. and follow the
multivariate normal distribution with mean zero and covariance matrix
$\mathbf{V}_N$.
The largest eigenvalue of the covariance matrix satisfies
$\| \mathbf{V}_N \| \le C $ for all $N\ge1$.
\end{enumerate}

Assumption (A1) deals with the properties of the sampling design. It
states that the sample size must be at least a positive fraction of the
population size,
that the one- and two-fold inclusion probabilities must be larger than
a positive number,
and that the two-fold inclusion probabilities should not be too far
from independence.
The latter is fulfilled, for example, for stratified sampling with
sampling without replacement within each stratum (Robinson and S\"
{a}rndal \cite{RoSa83}) and is discussed in details in H\`ajek \cite{Ha81} for
rejective sampling and other unequal probability sampling designs.
Assumption (A2) imposes H\"older continuity on the trajectories, a mild
regularity condition.
Assumption (A3) states that the design points have a quasi-uniform repartition
(this holds in particular for equidistant designs and designs generated
by a regular density function)
and that the grid size is essentially negligible compared to the
population size (e.g., if $d_N \propto N^{\alpha}$ for some
$\alpha\in(0,1)$).
In fact, the results of this paper also hold if $d_N/N$ stays bounded
away from zero and infinity as $N\to\infty$ (see Section \ref{sec5}).
Finally, (A4) imposes joint normality, short range temporal dependence,
and bounded variance for the measurement errors $\varepsilon_{kj}, 1\le j
\le d$.
It is trivially satisfied if the $\varepsilon_{kj} \sim N(0,\sigma_j^2)$
are independent with variances $\operatorname{Var} (\varepsilon_{kj}) \leq C$.
It is also verified if the $\varepsilon_{kj}$ arise from a discrete time
Gaussian process with short term temporal correlation such as ARMA or
stationary mixing processes.
Note that the Gaussian assumption is not central to our derivations: it
can be weakened and replaced by moment conditions on the error
distributions at the expense of much more complicated proofs.

%s3.1 ###
\subsection{Limit distribution of the Horvitz--Thompson
estimator}\label{AN HT estimator}

%Proceeding further, we would now like to
We now derive the asymptotic distribution of our estimator $\widehat
{\mu}_N$ in order to build asymptotic
confidence bands. Obtaining the asymptotic normality of estimators in
survey sampling is a technical and difficult issue even for simple
quantities such as means or totals of real numbers. Although confidence
intervals are commonly used in the survey sampling community, the
Central Limit Theorem (CLT) has only been checked rigorously, as far as
we know, for a few sampling designs. Erd\"os and R\'enyi \cite{ErRe59} and H\`
ajek \cite{Ha60} proved that the Horvitz--Thompson estimator is
asymptotically Gaussian for simple random sampling without replacement.
The CLT for rejective sampling is shown by H\`ajek \cite{Ha64} whereas the
CLT for other proportional to size sampling designs is studied by
Berger \cite{Be98}. Recently, these results were extended for some
particular cases of two-phase sampling designs (Chen and Rao \cite{ChRa07}).
Let us assume that the Horvitz--Thompson estimator satisfies a CLT for
real-valued quantities.
\begin{enumerate}[(A5)]
\item[(A5)] (\textit{Univariate CLT})
For any fixed $t\in[0,T]$,
%any bounded triangular array of fixed real numbers $\{ (a_{N,1}),
% letting $M_N = \sum_{k=1}^N a_{N,k}$ and $\widehat{ t}_{N}$ be the
%Horvitz--Thompson estimator of $M_N$ where
% the $n_N$ units are chosen randomly according to the sequence of
%fixed-size sampling designs $p_N$,
it holds that
\[
 \frac{ \widehat{\mu}_{N}(t) - \mu_N(t) }{\sqrt{
\operatorname{Var} ( \widehat{\mu}_{N}(t)  ) }}\leadsto N(0,1)
\]
as $N\to\infty$, where $\leadsto$ stands for convergence in distribution.
\end{enumerate}

We recall here the definition of the weak convergence in $ C([0,T])$
equipped with the supremum norm $\|\cdot\|_\infty$ (e.g., van der
Vaart and Wellner \cite{VaWe00}).
A sequence $(\xi_N)$ of random elements of $C([0,T])$ is said to
converge weakly to a limit $\xi$ in $ C([0,T]) $
if $\mathbb{E}(\phi(\xi_N)) \to\mathbb{E}(\phi(\xi))$ as $N\to
\infty$ for all bounded, uniformly continuous functionals $\phi$ on \mbox{$
(C([0,T]),\|\cdot\|_\infty)$}.

To establish the limit distribution of $ \widehat{\mu}_N$ in
$C([0,T])$, we need to assume the existence of a limit
% mean function
%$$ \mu(t) = \lim_{N\to\infty} \mu_N(t) $$ and
covariance function
\[
\gamma(s,t) = \lim_{N\to\infty} \frac{1}{N} \sum
_{k,l \in U_N} \Delta_{kl} \frac{X_k(s)}{\pi_k}
\frac{X_l(t)}{\pi_l}   .
\]
%
%{\textcolor{red}{Make sure that the convergence of $\gamma_N$ to $
%$C([0,T])$}}

In the following theorem, we state the asymptotic normality of the
estimator $\widehat{\mu}_N$ in the space $C([0,T])$ equipped with the
sup norm.

\begin{theo}\label{asymptotic normality mu_hat}
Assume \textup{(A1)--(A5)} and that %the bandwidth $h=h(N)$ satisfies
$\sqrt{N}h^\beta\to0$ and $dh/\log d \to\infty$ as $N\to\infty$.
Then
\[
\sqrt{N} ( \widehat{\mu}_N - \mu_N ) \leadsto G
\]
in $C([0,T])$, where %$\leadsto$ denotes weak convergence and
$G$ is a Gaussian process with mean zero and covariance function
$\gamma$.
\end{theo}

Theorem \ref{asymptotic normality mu_hat} provides a convenient way to
infer the local features of $\mu_N$. %Its proof is given in the
%Appendix.
It is applied in Section~\ref{SCB section} to the construction of
simultaneous confidence bands,
but it can also be used for a variety of statistical tests based on
supremum norms (see Degras \cite{De11}).

Observe that the conditions on the bandwidth $h$ and design size $d$
are not very constraining.
Suppose, for example, that $d \propto N^\eta$ and $h \propto N^{-\nu}$
for some $\eta, \nu>0$.
Then $d$ and $h$ satisfy the conditions of Theorem \ref{asymptotic
normality mu_hat} as soon as $(2 \beta)^{-1} < \nu< \eta< 1 $.
%Suppose that the design points are equispaced and $d_N = N^\eta,$ for
%some $\eta$ such that $(2 \beta)^{-1} <\eta<1$ with $\beta>1/2.$
%Then (A3) is automatically fulfilled and the conditions of Theorem~
Thus, for more regular trajectories, \textit{that is}, larger $\beta
,$ the bandwidth $h$ can be chosen with more flexibility.
%for choosing convenient values for $d$ and $h$.
%the number of design points and the size of the bandwidth.

%in that it allows to use statistics based on functionals of $\widehat{
%in $(C([0,T]),\|\cdot\|_\infty)$. This in turn allows to see whether
%the function $\mu_N$
%possesses a certain property everywhere on $[0,T]$, %e.g. positivity,
%monotonicity, or a parametric shape,
%as opposed to statistics based on euclidean distances that only check
%whether the property is satisfied globally
%(i.e. on average) on $[0,T]$. For instance, when testing for a
%parametric shape
%by measuring the distance between the fitted parametric and
%nonparametric curves,
%a supremum norm-based statistic can detect local discrepancies
%whereas a pseudo-likelihood ratio test or $L^2$-distance statistic
%is only sensitive to global agreement (see the Section 4 in).
The proof of Theorem \ref{asymptotic normality mu_hat} is similar in
spirit to that of Theorem 1 in Degras \cite{De11} and Proposition~3 in
Cardot and Josserand \cite{CaJo11}.\vadjust{\goodbreak}
Essentially, it breaks down into: (i) controlling uniformly on $[0,T]$
the bias of $\widehat{\mu}_N$, (ii) establishing the functional
asymptotic normality of the local linear smoother applied to the
sampled curves $X_k$ and (iii) controlling uniformly on $[0,T]$ (in
probability) the local linear smoother applied to the errors $\varepsilon_{jk}$. Part (i) is easily handled with standard results on
approximation properties of local polynomial estimators (see, e.g., Tsybakov \cite{Ts09}).
Part (ii) mainly consists in proving an asymptotic tightness property,
which entails the computation of entropy numbers and the use of maximal
inequalities (van der Vaart and Wellner \cite{VaWe00}). Part (iii) requires
first to show the finite-dimensional convergence of the smoothed error
process to zero and then to establish its tightness with similar
arguments as in part (ii). %The details of the proof are deferred to
%the Appendix.

%%%%%%%%%%%%%%%%%%%%%%%%%%%%%%%%%%%

%s3.2 ###
\subsection{Uniform consistency of the covariance estimator}\label{consistency in covariance estimation}

We first note that under (A1)--(A4), by the approximation properties of
local linear smoothers,
$\gamma_N$ converges uniformly to $\gamma$ on $[0,T]^2$ as $h \to0 $
and $N\to\infty$.
Hence, the consistency of $\widehat{\gamma}_N$ can be stated with
respect to $\gamma$ instead of $\gamma_N$.
In alignment with the related Proposition 2 in Cardot and Josserand
\cite{CaJo11} and Theorem 3 in Breidt and Opsomer \cite{BrOp00},
we need to make some assumption on the two-fold inclusion probabilities
of the sampling design $p_N$:
\begin{enumerate}[(A6)]
\item[(A6)]
\[\vspace*{-6pt}\lim_{N\to\infty} \max_{(k_1,k_2,k_3,k_4) \in
D_{4,N}} \big| \mathbb{E} \bigl\{
(I_{k_1}I_{k_2} - \pi_{k_1k_2}) (I_{k_3}I_{k_4}
- \pi_{k_3k_4}) \bigr\} \big| =0,\]
where $D_{4,N}$ is the set of all quadruples $(k_1,k_2,k_3,k_4)$ in
$U_N$ with distinct elements.
%of the 4-fold products $\mathbb{E} (\varepsilon_{k_1}\varepsilon_{k_2}
%where $(\varepsilon_1,\ldots,\varepsilon_d )\sim N(0,\mathbf{V}_N)$. NO! The
%normality of the errors already guarantees this condition.}}
\end{enumerate}
This assumption is discussed in detail in Breidt and Opsomer \cite{BrOp00} and
is fulfilled, for example, for stratified sampling.
%Also, to obtain uniform consistency, some additional rate condition on
%the bandwidth $h$ must be met.

\begin{theo}\label{teo2}
Assume \textup{(A1)--(A4)}, \textup{(A6)}, and that $h\to0$ and $dh^{1+\alpha} \to
\infty$ for some $\alpha>0$ as $N\to\infty$.
Then
%the covariance estimator $\widehat{\gamma}_N $ converges in
%probability uniformly to $\gamma$ over $[0,T]^2$, that is,
%
\[
\lim_{N\to\infty} \mathbb{E} \Bigl( \sup_{s,t\in[0,T]^2} \big| \widehat{
\gamma}_N (s,t)- \gamma(s,t) \big|^2 \Bigr) = 0, % \textrm{in probability.}
\]
where the expectation is jointly with respect to the design and the
multivariate normal model.
\end{theo}

  Note the additional condition on the bandwidth $h$ in Theorem \ref{teo2}.
If we suppose, as in the remark in Section~\ref{AN HT estimator},
that $d \propto N^\eta$ %$\eta$ such that
%$(2 \beta)^{-1} <\eta<1$ %with $\beta>1/2$
and $h \propto N^{-\nu}$ for some $(2 \beta)^{-1} < \nu< \eta< 1 $,
then condition $dh^{1+\alpha} \to\infty$ as $N\to\infty$ is
fulfilled with, for example, $\alpha= 1 - \eta/ 2 \nu$.
%as soon as $(2 \beta)^{-1} < \nu< \eta/(1+\alpha).$

%%%%%%%%%%%%%%%%%%%%%%%%%%%%%%%%%%%

%s3.3 ###
\subsection{Confidence bands}\label{SCB section}

In this section, we build confidence bands for $\mu_N$ of the form
\begin{equation}
\label{SCB} \biggl\{ \biggl[ \widehat{\mu}_N(t) \pm c
\frac{ \widehat{\sigma
}_N (t)}{N^{1/2}} \biggr],   t\in[0,T] \biggr\},
\end{equation}
where $c$ is a suitable number and $\widehat{\sigma}_N(t) = \widehat
{\gamma}_N(t,t)^{1/2}$.
More precisely, given a confidence level $1-\alpha\in(0,1),$ we seek
$c=c_\alpha$ that approximately satisfies
\begin{equation}
\label{SCB threshold} \mathbb{P} \bigl( \big|G(t)\big| \le c   \sigma(t) ,   \forall t
\in[0,T] \bigr) = 1-\alpha,
\end{equation}
where $G$ is a Gaussian process with mean zero and covariance function
$\gamma$, and where $\sigma(t) =\gamma(t,t)^{1/2}$.
Exact bounds for the supremum of Gaussian processes have been derived
for only a few particular cases (Adler and Taylor \cite{AdTa07}, Chapter 4).
Computing accurate and as explicit as possible bounds in a general
setting is a difficult issue and would require additional strong
conditions such as stationarity which have no reason to be fulfilled in
our setting.

% obtaining good explicit approximations is difficult and would require
%additional strong conditions such as stationarity
In view of Theorems \ref{asymptotic normality mu_hat}--\ref
{teo2} and Slutski's theorem,
the bands defined in \eqref{SCB} with $c$ chosen as in \eqref{SCB
threshold} will have approximate coverage level $1-\alpha$.
The following result provides a simulation-based method to compute $c$.

\begin{theo}\label{conditional weak functional convergence for
Gaussian processes}
Assume \textup{(A1)--(A6)} and %that $d^3h^5\to\infty$ as $N\to\infty$.
$dh^{1+\alpha} \to\infty$ for some $\alpha>0$ as $N\to\infty$.
Let $G$ be a Gaussian process with mean zero and covariance function
$\gamma$.
Let $(\widehat{G}_N)$ be a sequence of processes such that for each
$N$, conditionally on $\widehat{\gamma}_N$, $\widehat{G}_N$ is
Gaussian with mean zero and covariance $\widehat{\gamma}_N$ defined
in \eqref{covariance estimator}.
%Let also $(X_N)$ be a sequence of Gaussian processes with mean zero
%and respective covariance functions $\gamma_N$.
Then for all $c>0$, as $N\to\infty$, the following convergence holds
in probability:
\[
\mathbb{P} \bigl( \big|\widehat{G}_N(t)\big| \le c  \widehat{
\sigma}_N (t) ,   \forall t \in[0,T]   |   \widehat{
\gamma}_N \bigr) \to\mathbb{P} \bigl( \big|G(t)\big| \le c   \sigma(t),
\forall t \in [0,T] \bigr) .
\]
\end{theo}

  Theorem \ref{conditional weak functional convergence for
Gaussian processes} is derived by showing the weak convergence of
$(\widehat{G}_N)$ to $G$ in $C([0,T])$, which stems from Theorem \ref
{teo2} and the Gaussian nature of the
processes $\widehat{G}_N$. As in the first two theorems, maximal
inequalities are used to obtain the above weak convergence.
%fairly technical and appeals to involving random entropy numbers (see
%e.g. van der Vaart and Wellner (2000)).
%Its proof is postponed to the Appendix.
The practical importance of Theorem \ref{conditional weak functional
convergence for Gaussian processes} is that it allows to estimate the
number $c$ in \eqref{SCB threshold} via simulation (with the previous
notations): conditionally on $\widehat{\gamma}_N$, one can simulate a
large number of sample paths of the Gaussian process $ (\widehat{G}_N
/\widehat{\sigma}_N )$ and compute their supremum norms.
One then obtains a precise approximation to the distribution of $\|
\widehat{G}_N / \widehat{\sigma}_N \|_\infty$,
and it suffices to set $c$ as the quantile of order $(1-\alpha)$ of
this distribution:
\begin{equation}
\label{simulation threshold} \mathbb{P} \bigl( \big|\widehat{G}_N(t)\big|
\le c   \widehat{\sigma}_N (t) ,   \forall t \in[0,T]   |
\widehat{\gamma}_N \bigr) = 1-\alpha.
\end{equation}

\begin{coro}
Assume \textup{(A1)--(A6)}. Under the conditions of Theorems \ref{asymptotic
normality mu_hat}--\ref
{conditional weak functional convergence for Gaussian processes}, the
bands defined in \eqref{SCB} with the real $c=c(\widehat{\gamma}_N)$
chosen as in \eqref{simulation threshold} have asymptotic coverage
level $1-\alpha$, that is,
\[
\lim_{N\to\infty} \mathbb{P} \biggl( \mu_N(t)\in \biggl[ \widehat
{\mu}_N(t) \pm c   \frac{ \widehat{\sigma}_N (t)}{N^{1/2}} \biggr],  \forall t
\in[0,T] \biggr) = 1-\alpha.
\]
\end{coro}

%%%%%%%%%%%%%%%%%%%%%%%%%%%%%%%%%%%%%%%%%%%%%%%%%%%%%%%%%%%%%%%

%s4 ###
\section{A simulation study}\label{sec4}

In this section, we evaluate the performances of the mean curve
estimator as well as the coverage and the width of the confidence bands
for different bandwidth selection criteria and different levels of noise.
The simulations are conducted in the R environment.

%s4.1 ###
\subsection{Simulated data and sampling designs}
\begin{figure}[t]

\includegraphics{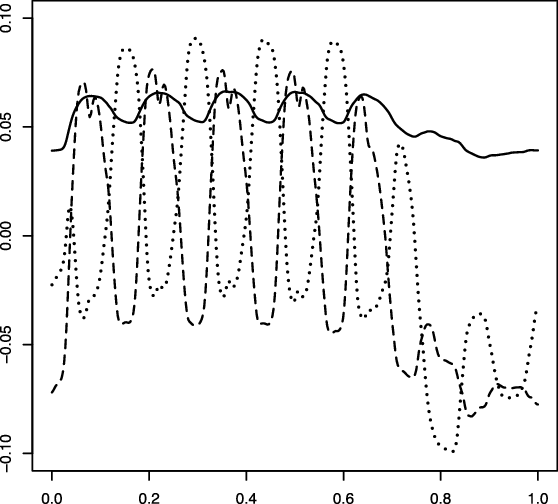}

\caption{Basis functions $v_1$ (solid line), $v_2$ (dashed line) and
$v_3$ (dotted line).}
\label{fig:eigen_vectors}
\end{figure}

\begin{figure}

\includegraphics{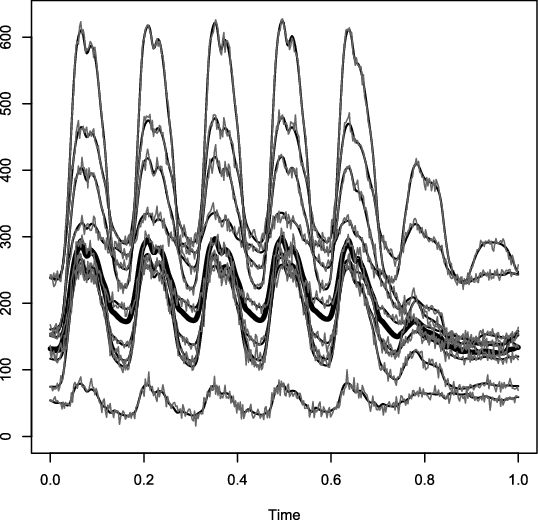}

\caption{A sample of $10$ curves for $\delta=0.05$ in the
heteroscedastic case. True trajectories are plotted with black lines
whereas noisy observations are plotted in gray. The mean profile is
plotted in bold line.}\vspace*{-3pt}
\label{fig:curveHET}
\end{figure}

We have generated a population of $N=20\mbox{,}000$ curves discretized at
$d=200$ and $d=400$ equidistant instants of time in $[0,1]$.
The curves of the population are generated so that they have
approximately the same distribution as the electricity consumption
curves analyzed in Cardot and Josserand \cite{CaJo11} and each individual
curve $X_k,$ for $k \in U,$ is simulated as follows
\begin{eqnarray}
X_k(t) = \mu(t) + \sum_{\ell=1}^3
Z_\ell v_\ell(t),\qquad    t \in[0,1],
\end{eqnarray}
where the mean function $\mu$ is drawn in Figure~\ref{fig:curveHET}
below
and the random variables $Z_\ell$ are independent realizations of a
centered Gaussian random variable with variance $\sigma_\ell^2$. The
three basis function $v_1, v_2$ and $v_3$ are orthonormal functions
which represent the main mode of variation of the signals, they are
represented in Figure \ref{fig:eigen_vectors}. Thus, the covariance
function of the population $\gamma(s,t)$ is simply
\begin{eqnarray}
\gamma(s,t) &=& \sum_{\ell=1}^3
\sigma_\ell^2  v_\ell(s) v_\ell(t).
\end{eqnarray}

%f1 ###

To select the samples, we have considered two probabilistic selection
procedures, with fixed sample size, $n=1000,$
\begin{itemize}
\item Simple random sampling without replacement (SRSWOR).
\item Stratified sampling with SRSWOR %simple random sampling without
%replacement
in all strata. The population $U$ is divided into a fixed number of
$H=5$ strata built by considering the quantiles $q_{0.5}, q_{0.7},
q_{0.85}$ and $q_{0.95}$ of the total consumption $\int_0^1 X_k(t) \,\mathrm{d}t$
for all units $k\in U$.
%of the population.
For example, the first strata contains all the units $k$ such that
$\int_0^1 X_k(t)\, \mathrm{d}t \leq q_{0.5},$ and thus its size is half of the
population size $N.$
The sample size $n_g$ in stratum $g$ is determined by a Neyman-like
allocation, as suggested in Cardot and Josserand \cite{CaJo11}, in order to
get a Horvitz--Thompson estimator of the mean trajectory whose variance
is as small as possible. The sizes of the different strata, which are
optimal according to this mean variance criterion, are reported in
Table \ref{tab:strata_description}.
\end{itemize}

%f2 ###

\begin{figure}%[t]

\includegraphics{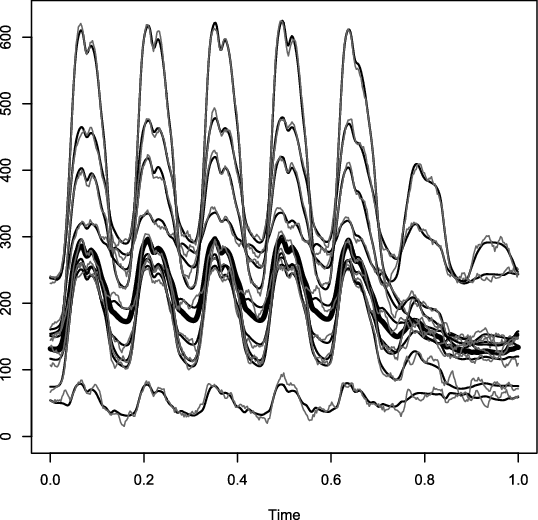}

\caption{A sample of $10$ curves for $\delta=0.05$ in the
autoregressive case. True trajectories are plotted with black lines
whereas noisy observations are plotted in gray. The mean profile is
plotted in bold line.}
\label{fig:curveAR3}
\end{figure}

\begin{figure}%[t]

\includegraphics{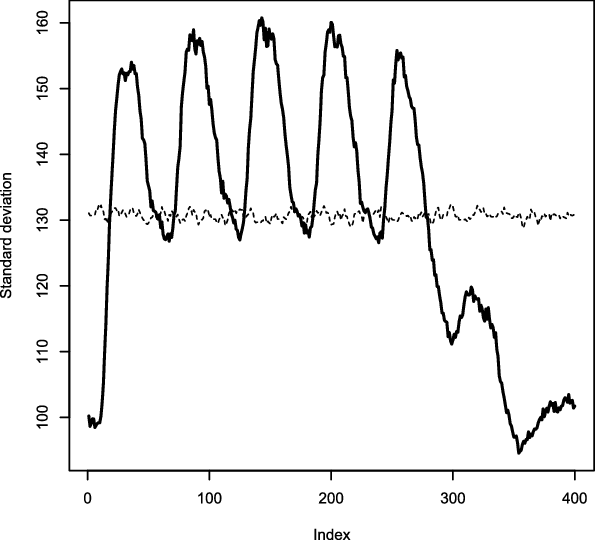}

\caption{Empirical standard deviation of the noise in the population
for $p=400$ discretization points. Standard deviation for
heteroscedastic case is drawn in solid line and dashed line for
correlated noise.}
\label{fig:sdnoise}
\end{figure}

We suppose we observe, for each unit $k$ in the sample $s,$ the
discretized trajectories, at $d$ equispaced points, $0=t_1 < \cdots<t_d=1,$
\begin{eqnarray}
\label{sim model} Y_{jk} = X_k(t_j) + \delta
\varepsilon_{jk}.
\end{eqnarray}
The parameter $\delta$ controls the noise level compared to the true
signal. We consider two different situations for the noise components
$\varepsilon_{jk}$:
\begin{itemize}
\item\textit{Heteroscedasticity}. The $\varepsilon_{jk} \sim N(0, \gamma
(t_j,t_j))$ are independent random variables
whose variances are proportional to the population variances at time $t_j$.
%with a variability which is directly related to the variability of the
%functions $X_k$ in the population.
%realizations of a centered Gaussian random variable with variance $$
%
\item\textit{Temporal dependence}. The $\varepsilon_{jk}$ are stationary
AR(3) processes with Gaussian innovations generated as follows
\[
\varepsilon_{jk} = 0.89 \varepsilon_{j-1,k} + 0.3
\varepsilon_{j-2,k} -0.4 \varepsilon_{j-3,k} + \eta_{jk}.
\]
The $\eta_{jk} \sim N(0, \sigma_\eta^2)$ are i.i.d. and $\sigma_\eta^2$ is such that $\mathbb{E}(\varepsilon_{jk}^2) = d^{-1} \sum_{j=1}^d
\gamma(t_j,t_j).$
% The $\eta_{jk} $ are {i.i.d.} realizations of a $N(0, \sigma_\eta^2)$
%where $\sigma_\eta^2$ is such that $\E(\varepsilon_{jk}^2) = d^{-1}
\end{itemize}
As an illustrative example, a sample of $n=10$ noisy discretized curves
are plotted in Figure~\ref{fig:curveHET} with heteroscedastic noise
components and in Figure~\ref{fig:curveAR3} for correlated noise. It
should be noted that the observed trajectories corrupted by the
correlated noise are much smoother %, with this particular
%autoregressive process for the noise,
than the trajectories corrupted by the heteroscedastic noise.
The empirical standard deviation in the population, for these two
different type of noise are drawn in Figure~\ref{fig:sdnoise}.

%t1 ###
\begin{table}[b]\vspace*{-3pt}
\caption{Strata sizes and optimal allocations}
\label{tab:strata_description}
\begin{tabular*}{\textwidth}{@{\extracolsep{\fill}}llllll@{}}
\hline
&\multicolumn{5}{l@{}}{Stratum number}\\[-5pt]
&\multicolumn{5}{l@{}}{\hrulefill}\\
 & 1 & 2 & 3 & 4 & 5 \\
\hline
Stratum size & 10,000 & 4000 & 3000 & 2000 & 1000 \\
Allocation & \phantom{10,}655 & \phantom{0}132 & \phantom{00}98 & \phantom{00}68 & \phantom{00}47 \\
\hline
\end{tabular*}
\end{table}

%s4.2 ###
\subsection{Weighted cross-validation for bandwidth selection}\label{sec4.2}

Assuming we can access the exact trajectories $X_k,   k\in s$ (which
is the case in simulations),
%$\mathbf{X}_k = (X_k(t_1), \ldots, X_k(t_d)),   k \in s$,
we consider the oracle-type estimator % that will be used as a
%benchmark:
%
\begin{equation}
\label{oracle} \widehat{\mu}_s = \sum_{k \in s}
\frac{X_k}{\pi_k}   ,
\end{equation}
which will be a benchmark in our numerical study.
We compare different interpolation and smoothing strategies for
estimating the $X_k,   k\in s$:
\begin{itemize}
\item Linear interpolation of the $Y_{jk}$ as in Cardot and Josserand \cite{CaJo11}.
\item Local linear smoothing of the $Y_{jk}$ with bandwidth $h$ as in
(\ref{def:smoother}).
% in order to get the functional estimations $\widehat{X}_k(t), t \in
%[0,1]$ and then compute $\widehat{\mu}$ as in (\ref{def:smoothmean}).
\end{itemize}

  The crucial parameter here is $h$. To evaluate the interest
of smoothing and the performances of data-driven bandwidth selection criteria,
we consider an error measure that
%also consider another benchmark which
compares the oracle $\widehat{\mu}_s$ to any estimator $\widehat{\mu
}$ based on the noisy data $Y_{jk},   k \in s,   j=1, \ldots, d$:
\begin{equation}
\label{eq_Lh} L(\widehat{\mu}) = \int_0^T
\bigl(\widehat{\mu}(t) - \widehat {\mu}_s(t) \bigr)^2\, \mathrm{d}t
  .
\end{equation}
%
% \begin{equation}
% \label{eq_Lh}
% L(\widehat{\mu}) &= \frac{1}{d} \left\| \widehat{\bolds{\mu}}_s
%- \widehat{\bolds{\mu}} \right\|^2,
% \end{equation}
% where $\widehat{\bolds{\mu}} = (\widehat{\mu}(t_1), \ldots,
%norm in $\mathbb{R}^d.$
Considering the estimator defined in (\ref{def:smoothmean}), we denote
by $h_{\mathrm{oracle}}$ the bandwidth $h$ that minimizes \eqref
{eq_Lh}. The mean estimator built with bandwidth $h_{\mathrm{oracle}}$
is called smooth oracle estimator. %the corresponding estimator.
%Note that in practice, we cannot carry out this estimator, called, but
%it will be the best possible estimator of $\widehat{\mu}_s$ given the
%data.

%f3 ###

%f4 ###

When $\sum_{k \in s} \pi_k^{-1} = N$, as in SRSWOR and stratified sampling,
it can be easily checked that $\widehat{\mu}_s$
is the minimum argument of the weighted least squares functional
\begin{equation}\label{def:wopest}
\sum_{k \in s} w_k \int
_0^T \bigl(X_k(t) - \mu(t)
\bigr)^2 \,\mathrm{d}t
\end{equation}
with respect to $\mu\in L^2([0,T]),$ where the weights are $w_k =
(N\pi_k)^{-1}.$
%Let us define $\mathbf{Y}_k = (Y_k(t_1), \ldots, Y_k(t_d))$ and $
Then, a simple and natural way to select bandwidth $h$ is to consider
the following design-based cross validation
\begin{equation}
\label{def:WCV} \operatorname{WCV}(h) = \sum_{k \in s}
w_k\sum_{j=1}^d \bigl(
Y_{jk} - \widehat{\mu}_N^{-k}(t_j)
\bigr)^2 ,
\end{equation}
where
\[
\widehat{\mu}_N^{-k} (t)= \sum
_{\ell\in s, \ell\neq k} \widetilde {w}_{\ell k}  \widehat{X}_\ell(t),
\]
with new weights $\widetilde{w}_{\ell k}.$ A heuristic justification
for this approach is that, given $s,$ we have $\mathbb{E} [
\varepsilon_{jk}(X_k(t_j) - \widehat{\mu}_N^{-k}(t_j)) | s ]=0$ for
$j=1,\ldots, d$ and $k \in s.$ Thus,
\begin{eqnarray*}
\mathbb{E} \bigl[ \operatorname{WCV}(h) | s \bigr] &=& \sum_{k\in s}
w_k \sum_{j=1}^d \bigl\{
\mathbb{E} \bigl[ \bigl(X_k(t_j) - \widehat{\mu
}_N^{-k}(t_j) \bigr)^2 | s \bigr]\\
&&{}
+ 2\mathbb{E} \bigl[ \varepsilon_{jk}\bigl(X_k(t_j)
- \widehat{\mu}_N^{-k}(t_j)\bigr) | s \bigr]
+ \mathbb {E} \bigl[ \varepsilon_{jk}^2 \bigr] \bigr\}
\\
&=& \sum_{k \in s} w_k \sum
_{j=1}^d \mathbb{E} \bigl[ \bigl(X_k(t_j)
- \widehat {\mu}_N^{-k}(t_j)
\bigr)^2 | s \bigr] + \operatorname{tr}(\mathbf{V}_N) % &= \sum_{k \in s} w_k \ \left\| \mathbf{X}_k - \widehat{\boldsymbol{
\end{eqnarray*}
and, up to $ \operatorname{tr}(\mathbf{V}_N)$ which does not depend on $h$,
the minimum value of the expected cross validation criterion should be
attained for estimators which are not too far from $\widehat{\mu}_s$.
%with the notations of \eqref{def tilde X eps Delta} and $\tilde{
%which, up to constant term $\operatorname{tr}(\mathbf{V}_n),$ is very similar
%to (\ref{def:wopest}).

This weighted cross validation criterion is simpler than the cross
validation criteria based on the estimated variance proposed in Opsomer
and Miller \cite{OpMi05}.
% or Berger and Skinner (2005).
Indeed, in our case, the bias may be non-negligible and focusing only
on the variance part of the error leads to too large selected values
for the bandwidth. Furthermore, Opsomer and Miller \cite{OpMi05} suggested to
consider weights defined as follows $\widetilde{w}_{\ell k} = w_\ell
/(1 - w_k)$. For SRSWOR, since $w_k = n^{-1}$ one has $\widetilde
{w}_{\ell k} = (n-1)^{-1},$ so that the weighted cross validation
criterion defined in (\ref{def:WCV}) is exactly the cross validation
criterion introduced by Rice and Silverman \cite{RiSi91} in the independent
case. We denote in the following by $h_{\mathrm{cv}}$ the bandwidth
value minimizing this criterion.

For stratified sampling, a better approximation which keeps the
design-based properties of the estimator $\widehat{\mu}_N^{-k}$ can
be obtained by taking into account the sampling rates in the different
strata. Assume the population $U$ is partitioned in strata $U_\nu$
of respective sizes $N_\nu,\nu=1, \ldots, H,$ and we sample $n_\nu$
observations in each $U_\nu$ by SRSWOR. If $k\in U_\nu,$ we have $w_k
= N_\nu(N n_\nu)^{-1}.$ Thus, we take $\widetilde{w}_{\ell k} =
(N_\nu-1)\{(N-1) (n_\nu-1)\}^{-1}$ for all the units $\ell\neq k$ in
stratum $U_\nu$ and scale the weights for all the units $\ell'$ of
the sample that do not belong to stratum $g,$ $\widetilde{w}_{\ell'
k}= N(N-1)^{-1} w_{\ell'}$. We denote by $h_{\mathrm{wcv}}$ the
bandwidth value minimizing (\ref{def:WCV}).

%s4.3 ###
\subsection{Estimation errors and confidence bands}

We draw $1000$ samples in the population of curves and compare the
different estimators of Section \ref{sec4.2}
with the $L^2$ loss criterion
\begin{equation}
R(\widehat{\mu}) = \int_0^T \bigl( \widehat{
\mu}(t) - \mu(t) \bigr)^2 \,\mathrm{d}t
\end{equation}
for different values of $\delta$ and $d$ in \eqref{sim model}. For
comparison, the quadratic approximation error for function $\mu$ by
its average value, $\overline{\mu} = T^{-1}\int_0^T \mu(t) \,\mathrm{d}t,$ is
$R(\overline{\mu}) = 3100.$

%t2 ###
\begin{table}%[t]
\caption{(Heteroscedastic noise). Estimation errors according to
$R(\widehat{\mu})$ for different noise levels and bandwidth choices,
with $d=200$ observation times. Units are selected by SRSWOR or
stratified sampling}
\label{tab:summary_mu_200}
\begin{tabular*}{\textwidth}{@{\extracolsep{\fill}}llllllllll@{}}
\hline
& & \multicolumn{4}{l}{SRSWOR} & \multicolumn{4}{l@{}}{Stratified sampling} \\[-5pt]
& & \multicolumn{4}{l}{\hrulefill} & \multicolumn{4}{l@{}}{\hrulefill} \\
$\delta$ & $h$ & Mean & 1Q & Median & 3Q & Mean & 1Q & Median & 3Q \\
\hline
5\% & lin & 17.65 & 3.08 & 8.73 & 23.50 & 4.22 & 1.44 & 2.79 & 5.59 \\
& $h_\mathrm{cv}$ & 17.65 & 3.07 & 8.71 & 23.51 & 6.49 & 3.61 & 5.36 &8.03 \\
& $h_\mathrm{wcv}$ & 17.65 & 3.07 & 8.71 & 23.51 & 4.22 & 1.45 & 2.78& 5.56 \\
& $h_\mathrm{oracle}$ & 17.65 & 3.07 & 8.72 & 23.50 & 4.22 & 1.45 &2.78 & 5.57 \\
& $\widehat{\mu}_s$ & 17.60 & 3.01 & 8.70 & 23.36 & 4.17 & 1.38 &2.76 & 5.55 \\
[6pt]
25\% & lin & 17.69 & 3.94 & 8.99 & 21.52 & 5.26 & 2.63 & 4.15 & 6.54 \\
& $h_\mathrm{cv}$ & 17.53 & 3.83 & 8.76 & 21.53 & 6.98 & 4.29 & 5.83 &8.47 \\
& $h_\mathrm{wcv}$ & 17.53 & 3.83 & 8.76 & 21.53 & 5.02 & 2.39 & 3.89& 6.33 \\
& $h_\mathrm{oracle}$ & 17.52 & 3.81 & 8.78 & 21.52 & 5.01 & 2.37 &3.88 & 6.27 \\
& $\widehat{\mu}_s$ & 16.58 & 2.85 & 7.87 & 20.01 & 4.07 & 1.46 &2.94 & 5.28 \\
\hline
\end{tabular*}
\end{table}

% \begin{table}[t]
% \centering
% \small
% \begin{tabular}{|c|l|c|c|c|c||c|c|c|c|}
% \hline
% \multicolumn{2}{|c|}{} & \multicolumn{4}{c||}{SRSWOR} &
% \hline
% $\delta$ & $h$ & Mean & 1Q & Median & 3Q & Mean & 1Q & Median & 3Q \\
% 5\% & lin & 16.50 & 3.12 & 8.55 & 21.40 & 4.09 & 1.43 & 2.84 & 5.19 \\
%& $h_\mathrm{cv}$ & 17.95 & 4.61 & 10.26 & 22.93 & 6.02 & 3.24 & 4.80
%& 7.50 \\
%& $h_\mathrm{wcv}$ & 17.95 & 4.61 & 10.26 & 22.93 & 4.33 & 1.68 & 3.09
%& 5.46 \\
%& $h_\mathrm{oracle}$ & 17.95 & 4.61 & 10.26 & 22.93 & 4.33 & 1.68 &
%3.09 & 5.46 \\
%& $\widehat{\mu}_s$ & 16.46 & 3.11 & 8.47 & 21.55 & 4.04 & 1.41 & 2.78
%& 5.08 \\
%25\% & lin & 17.40 & 4.13 & 9.24 & 22.25 & 5.19 & 2.63 & 3.95 & 6.57 \\
% & $h_\mathrm{cv}$ & 17.84 & 4.53 & 9.76 & 22.63 & 7.13 & 4.46 & 6.05
%& 8.61 \\
% & $h_\mathrm{wcv}$ & 17.84 & 4.53 & 9.76 & 22.63& 5.39 & 2.76 & 4.24
%& 6.74 \\
% & $h_\mathrm{oracle}$ & 17.84 & 4.53 & 9.76 & 22.63 & 5.39 & 2.76 &
%4.24 & 6.74 \\
% & $\widehat{\mu}_s$ & 16.39 & 3.13 & 8.23 & 21.54 & 3.96 & 1.49 &
%2.87 & 5.19 \\
% \hline
% \end{tabular}
% \caption{(correlated noise). Estimation errors according to $R(
%$d=200$ %observation times. Units are selected by SRSWOR or stratified
%sampling.}
% \label{tab:summary_mu_200c}
% \end{table}

%t3 ###
\begin{table}[b]
\caption{(Heteroscedastic noise). Estimation errors according to
$R(\widehat{\mu})$ for different noise levels and bandwidth choices,
with $d=400$ observation times. Units are selected by SRSWOR or
stratified sampling}
\label{tab:summary_mu}
\begin{tabular*}{\textwidth}{@{\extracolsep{\fill}}llllllllll@{}}
\hline
& & \multicolumn{4}{l}{SRSWOR} & \multicolumn{4}{l@{}}{Stratified sampling} \\[-5pt]
& & \multicolumn{4}{l}{\hrulefill} & \multicolumn{4}{l@{}}{\hrulefill} \\
$\delta$ & $h$ & Mean & 1Q & Median & 3Q & Mean & 1Q & Median & 3Q \\
\hline
5\% & lin & 18.03 & 3.39 & 9.24 & 23.27 & 4.05 & 1.45 & 2.86 & 5.35 \\
& $h_\mathrm{cv}$ & 18.02 & 3.38 & 9.26 & 23.34 & 6.09 & 3.24 & 4.87 &7.56 \\
& $h_\mathrm{wcv}$ & 18.02 & 3.38 & 9.26 & 23.34 & 4.05 & 1.45 & 2.82& 5.40 \\
& $h_\mathrm{oracle}$ & 18.02 & 3.38 & 9.27 & 23.32 & 4.04 & 1.43 &2.83 & 5.39 \\
& $\widehat{\mu}_s$ & 17.98 & 3.35 & 9.20 & 23.17 & 4.00 & 1.39 &2.81 & 5.29 \\
[6pt]
25\% & lin & 18.16 & 3.89 & 9.43 & 22.86 & 5.25 & 2.85 & 4.24 & 6.57 \\
& $h_\mathrm{cv}$ & 17.55 & 3.30 & 8.89 & 22.09 & 6.45 & 3.77 & 5.37 &8.11 \\
& $h_\mathrm{wcv}$ & 17.55 & 3.30 & 8.89 & 22.09 & 4.57 & 2.12 & 3.49& 5.81 \\
& $h_\mathrm{oracle}$ & 17.55 & 3.28 & 8.89 & 22.09 & 4.56 & 2.11 &3.48 & 5.81 \\
& $\widehat{\mu}_s$ & 17.04 & 2.75 & 8.38 & 21.87 & 4.04 & 1.60 &3.02 & 5.31 \\
\hline
\end{tabular*}
\end{table}

%t4 ###
\begin{table}%[t]
\caption{(Correlated noise). Estimation errors according to
$R(\widehat{\mu})$ for different noise levels and bandwidth choices,
with $d=400$ observation times. Units are selected by SRSWOR or
stratified sampling}
\label{tab:summary_mu_200c}
\begin{tabular*}{\textwidth}{@{\extracolsep{\fill}}llllllllll@{}}
\hline
& & \multicolumn{4}{l}{SRSWOR} & \multicolumn{4}{l@{}}{Stratified sampling} \\[-5pt]
& & \multicolumn{4}{l}{\hrulefill} & \multicolumn{4}{l@{}}{\hrulefill} \\
$\delta$ & $h$ & Mean & 1Q & Median & 3Q & Mean & 1Q & Median & 3Q \\
\hline
5\% & lin & 16.23 & 3.05 & 8.67 & 20.86 & 4.08 & 1.40 & 2.88 & 5.44 \\
& $h_\mathrm{cv}$ & 16.24 & 3.07 & 8.66 & 20.88 & 5.90 & 2.99 & 4.70 &7.33 \\
& $h_\mathrm{wcv}$ & 16.24 & 3.07 & 8.66 & 20.88 & 4.10 & 1.38 & 2.90& 5.47 \\
& $h_\mathrm{oracle}$ & 16.24 & 3.06 & 8.65 & 20.88 & 4.10 & 1.38 &2.90 & 5.46 \\
& $\widehat{\mu}_s$ & 16.19 & 3.01 & 8.69 & 20.86 & 4.04 & 1.34 &2.82 & 5.36 \\
[6pt]
25\% & lin & 17.18 & 3.88 & 9.38 & 22.04 & 5.22 & 2.65 & 4.07 & 6.47 \\
& $h_\mathrm{cv}$ & 17.13 & 3.84 & 9.28 & 22.02 & 6.76 & 3.98 & 5.76 &8.32 \\
& $h_\mathrm{wcv}$ & 17.13 & 3.84 & 9.28 & 22.02 & 5.16 & 2.59 & 4.02& 6.37 \\
& $h_\mathrm{oracle}$ & 17.12 & 3.81 & 9.25 & 22.02 & 5.15 & 2.59 &4.01 & 6.37 \\
& $\widehat{\mu}_s$ & 16.12 & 2.87 & 8.17 & 21.00 & 4.04 & 1.49 &2.94 & 5.27 \\
\hline
\end{tabular*}
\end{table}

The empirical mean as well as the first, second and third quartiles of
the estimation error $R(\widehat{\mu})$ are given, when $d=200$, in
Table~\ref{tab:summary_mu_200} for the heteroscedastic noise case.
Results for $d=400$ are presented for the heteroscedastic case in
Table~\ref{tab:summary_mu} and in Table~\ref{tab:summary_mu_200c} for
the correlated case.

We first note that in all simulations, stratified sampling largely
improves the estimation of the mean curve in comparison to SRSWOR.
Also, linear interpolation performs nearly as well as the smooth oracle
estimator for large samples,
especially when the noise level is low ($\delta= 5\% $).
As far as bandwidth selection is concerned, the usual cross validation
criterion $h_{\mathrm{cv}}$
is not adapted to unequal probability sampling and tends to select too
large bandwidth values.
In particular, it does not perform as well as linear interpolation for
stratified sampling.
On the other hand, our weighted cross-validation method seems effective
for selecting the bandwidth.
It produces estimators that are very close to the oracle and that dominate
the other estimators when the noise level is moderate or high ($\delta
= 25\% $).

%t5 ###
\begin{table}[b]
\caption{(Heteroscedastic noise). Estimation errors according to
$L(\widehat{\mu})$ for different noise levels and bandwidth choices,
with $d=200$ observation times. Units are selected by SRSWOR or
stratified sampling}\label{tab:summary_Lh_200}
\begin{tabular*}{\textwidth}{@{\extracolsep{\fill}}llllllllll@{}}
\hline
& & \multicolumn{4}{l}{SRSWOR} & \multicolumn{4}{l@{}}{Stratified sampling} \\[-5pt]
& & \multicolumn{4}{l}{\hrulefill} & \multicolumn{4}{l@{}}{\hrulefill} \\
$\delta$ & $h$ & Mean & 1Q & Median & 3Q & Mean & 1Q & Median & 3Q \\
\hline
5\% & lin & 0.044 & 0.041 & 0.044 & 0.047 & 0.049 & 0.046 & 0.049 &0.053 \\
& $h_\mathrm{cv}$ & 0.044 & 0.041 & 0.044 & 0.048 & 2.520 & 2.083 &2.852 & 3.032 \\
& $h_\mathrm{wcv}$ & 0.044 & 0.041 & 0.044 & 0.048 & 0.058 & 0.054 &0.058 & 0.062 \\
& $h_\mathrm{oracle}$ & 0.044 & 0.041 & 0.044 & 0.047 & 0.049 & 0.045& 0.049 & 0.052 \\
[6pt]
25\% & lin & 1.087 & 1.011 & 1.080 & 1.156 & 1.214 & 1.134 & 1.210 &1.287 \\
& $h_\mathrm{cv}$ & 0.905 & 0.837 & 0.901 & 0.970 & 3.155 & 2.638 &3.260 & 3.602 \\
& $h_\mathrm{wcv}$ & 0.905 & 0.837 & 0.901 & 0.970 & 1.009 & 0.936 &1.004 & 1.076 \\
& $h_\mathrm{oracle}$ & 0.898 & 0.830 & 0.894 & 0.962 & 0.990 & 0.919& 0.988 & 1.055 \\
\hline
\end{tabular*}
\end{table}

This is clearer when we look at criterion $L(\widehat{\mu}),$ defined
in (\ref{eq_Lh}), which only focuses on the part of the estimation
error which is due to the noise.
Results are presented in Table~\ref{tab:summary_Lh_200} for $d=200$ in
the heteroscedastic case. For $d=400,$ errors are given in Table~\ref
{tab:summary_Lh} in the heteroscedastic case and in Table~\ref
{tab:summary_Lhc} for correlated noise. When the noise level is high,
we observe a significant impact of the number of discretization points
on the accuracy of the smoothed estimators. Our individual
trajectories, which have roughly the same shape as load curves, are
actually not very smooth so that smoothing approaches are only really
interesting, compared to linear interpolation, when the number of
discretization points $d$ is large enough. Finally, it also becomes
clearer that a key parameter is the bandwidth value which has to be
chosen with appropriate criteria that must take the sampling weights
into account. When the noise level is low ($\delta=5\%$), the error
according to criterion $L(\widehat{\mu})$ is multiplied by at least
15 in stratified sampling.

%t6 ###
\begin{table}%[t]
\caption{(Heteroscedastic noise). Estimation errors according to
$L(\widehat{\mu})$ for different noise levels and bandwidth choices,
with $d=400$ observation times. Units are selected by SRSWOR or
stratified sampling}
\label{tab:summary_Lh}
\begin{tabular*}{\textwidth}{@{\extracolsep{\fill}}llllllllll@{}}
\hline
& & \multicolumn{4}{l}{SRSWOR} & \multicolumn{4}{l@{}}{Stratified sampling} \\[-5pt]
& & \multicolumn{4}{l}{\hrulefill} & \multicolumn{4}{l@{}}{\hrulefill} \\
$\delta$ & $h$ & Mean & 1Q & Median & 3Q & Mean & 1Q & Median & 3Q \\
\hline
5\% & lin & 0.044 & 0.042 & 0.044 & 0.047 & 0.049 & 0.047 & 0.049 &0.051 \\
& $h_\mathrm{cv}$ & 0.040 & 0.038 & 0.040 & 0.042 & 2.231 & 1.612 &1.917 & 2.806 \\
& $h_\mathrm{wcv}$ & 0.040 & 0.038 & 0.040 & 0.042 & 0.052 & 0.049 &0.052 & 0.055 \\
& $h_\mathrm{oracle}$ & 0.040 & 0.038 & 0.040 & 0.042 & 0.044 & 0.041& 0.044 & 0.046 \\
[6pt]
25\% & lin & 1.089 & 1.030 & 1.087 & 1.142 & 1.219 & 1.155 & 1.212 &1.280 \\
& $h_\mathrm{cv}$ & 0.498 & 0.462 & 0.495 & 0.535 & 2.591 & 1.932 &2.344 & 3.254 \\
& $h_\mathrm{wcv}$ & 0.498 & 0.462 & 0.495 & 0.535 & 0.552 & 0.509 &0.549 & 0.594 \\
& $h_\mathrm{oracle}$ & 0.497 & 0.460 & 0.494 & 0.533 & 0.547 & 0.505& 0.545 & 0.586 \\
\hline
\end{tabular*}
\end{table}

%t7 ###
\begin{table}[b]
\caption{(Correlated noise). Estimation errors according to
$L(\widehat{\mu})$ for different noise levels and bandwidth choices,
with $d=400$ observation times. Units are selected by SRSWOR or
stratified sampling}
\label{tab:summary_Lhc}
\begin{tabular*}{\textwidth}{@{\extracolsep{\fill}}llllllllll@{}}
\hline
& & \multicolumn{4}{l}{SRSWOR} & \multicolumn{4}{l@{}}{Stratified sampling} \\[-5pt]
& & \multicolumn{4}{l}{\hrulefill} & \multicolumn{4}{l@{}}{\hrulefill} \\
$\delta$ & $h$ & Mean & 1Q & Median & 3Q & Mean & 1Q & Median & 3Q \\
\hline
5\% & lin & 0.17 & 0.15 & 0.16 & 0.18 & 0.05 & 0.04 & 0.04 & 0.05 \\
& $h_\mathrm{cv}$ & 0.17 & 0.15 & 0.16 & 0.18 & 1.94 & 1.53 & 1.59 &2.90 \\
& $h_\mathrm{wcv}$ & 0.17 & 0.15 & 0.16 & 0.18 & 0.07 & 0.07 & 0.07 &0.08 \\
& $h_\mathrm{oracle}$ & 0.17 & 0.15 & 0.16 & 0.18 & 0.07 & 0.07 & 0.07& 0.08 \\
[6pt]
25\% & lin & 1.09 & 1.03 & 1.09 & 1.14 & 1.20 & 1.08 & 1.19 & 1.32 \\
& $h_\mathrm{cv}$ & 0.50 & 0.46 & 0.50 & 0.53 & 2.83 & 2.19 & 2.57 &3.67 \\
& $h_\mathrm{wcv}$ & 0.50 & 0.46 & 0.50 & 0.53 & 1.15 & 1.02 & 1.13 &1.26 \\
& $h_\mathrm{oracle}$ & 0.49 & 0.46 & 0.49 & 0.53 & 1.13 & 1.01 & 1.12& 1.25 \\
\hline
\end{tabular*}
\end{table}

%t8 ###
\begin{table}%[t]
\tabcolsep=4pt
\caption{(Heteroscedastic noise). Empirical coverage levels
$1-\widehat{\alpha}$ and confidence band areas for different noise
levels and bandwidth choices, with $d=200$ observation times. Units are
selected by SRSWOR or stratified sampling}
\label{tab:confidence_bands_200}
\begin{tabular*}{\textwidth}{@{\extracolsep{\fill}}llllllllllll@{}}
\hline
& & \multicolumn{5}{l}{SRSWOR} & \multicolumn{5}{l@{}}{Stratified sampling} \\[-5pt]
& & \multicolumn{5}{l}{\hrulefill} & \multicolumn{5}{l@{}}{\hrulefill} \\
$\delta$ & $h$ & $1-\widehat{\alpha}$ & Mean & 1Q & Median & 3Q &$1-\widehat{\alpha}$ & Mean & 1Q & Median & 3Q\\
\hline
5\% & lin & 97.2 & 10.91 & 10.74 & 10.90 & 11.07 & 98.1 & 5.95 & 5.87 &5.95 & 6.02 \\
& $h_\mathrm{cv}$ & 97.3 & 10.89 & 10.73 & 10.89 & 11.06 & 47.5 & 5.68& 5.60 & 5.68 & 5.76 \\
& $h_\mathrm{wcv}$ & 97.3 & 10.89 & 10.73 & 10.89 & 11.06 & 97.5 &5.92 & 5.84 & 5.91 & 6.00 \\
& $h_\mathrm{oracle}$ & 97.2 & 10.90 & 10.72 & 10.90 & 11.07 & 98.0 &5.94 & 5.86 & 5.94 & 6.02 \\
& $\widehat{\mu}_s$ & 97.3 & 10.54 & 10.36 & 10.54 & 10.70 & 98.2 &5.59 & 5.51 & 5.60 & 5.67 \\
[6pt]
25\% & lin & 97.7 & 13.23 & 13.06 & 13.22 & 13.41 & 98.3 & 8.27 & 8.19& 8.27 & 8.36 \\
& $h_\mathrm{cv}$ & 97.2 & 12.66 & 12.49 & 12.65 & 12.83 & 64.7 & 6.70& 6.60 & 6.69 & 6.79 \\
& $h_\mathrm{wcv}$ & 97.2 & 12.66 & 12.49 & 12.65 & 12.83 & 97.3 &7.56 & 7.48 & 7.56 & 7.65 \\
& $h_\mathrm{oracle}$ & 97.3 & 12.70 & 12.50 & 12.70 & 12.87 & 97.5 &7.68 & 7.58 & 7.68 & 7.79 \\
& $\widehat{\mu}_s$ & 97.0 & 10.53 & 10.37 & 10.52 & 10.70 & 97.7 &5.59 & 5.51 & 5.59 & 5.66 \\
\hline
\end{tabular*}
\end{table}

%t9 ###
\begin{table}[b]
\tabcolsep=4pt
\caption{(Heteroscedastic noise). Empirical coverage levels
$1-\widehat{\alpha}$ and confidence band areas for different noise
levels and bandwidth choices, with $d=400$ observation times. Units are
selected by SRSWOR or stratified sampling}
\label{tab:confidence_bands}
\begin{tabular*}{\textwidth}{@{\extracolsep{\fill}}llllllllllll@{}}
\hline
& & \multicolumn{5}{l}{SRSWOR} & \multicolumn{5}{l@{}}{Stratified sampling} \\[-5pt]
& & \multicolumn{5}{l}{\hrulefill} & \multicolumn{5}{l@{}}{\hrulefill} \\
$\delta$ & $h$ & $1-\widehat{\alpha}$ & Mean & 1Q & Median & 3Q &$1-\widehat{\alpha}$ & Mean & 1Q & Median & 3Q\\
\hline
5\% & lin & 97.7 & 10.79 & 10.63 & 10.79 & 10.95 & 97.9 & 6.03 & 5.95 &6.02 & 6.11 \\
& $h_\mathrm{cv}$ & 97.6 & 10.76 & 10.59 & 10.77 & 10.92 & 48.4 & 5.64& 5.57 & 5.63 & 5.72 \\
& $h_\mathrm{wcv}$ & 97.6 & 10.76 & 10.59 & 10.77 & 10.92 & 97.6 &5.89 & 5.82 & 5.89 & 5.97 \\
& $h_\mathrm{oracle}$ & 97.6 & 10.76 & 10.59 & 10.77 & 10.92 & 97.6 &5.96 & 5.88 & 5.96 & 6.04 \\
& $\widehat{\mu}_s$ & 97.7 & 10.50 & 10.33 & 10.50 & 10.65 & 97.8 &5.60 & 5.52 & 5.59 & 5.68 \\
[6pt]
25\% & lin & 97.6 & 12.69 & 12.52 & 12.70 & 12.86 & 98.3 & 8.59 & 8.49& 8.59 & 8.68 \\
& $h_\mathrm{cv}$ & 97.5 & 12.47 & 12.31 & 12.48 & 12.64 & 58.1 & 6.34& 6.24 & 6.34 & 6.44 \\
& $h_\mathrm{wcv}$ & 97.5 & 12.47 & 12.31 & 12.48 & 12.64 & 97.6 &7.09 & 7.00 & 7.08 & 7.17 \\
& $h_\mathrm{oracle}$ & 97.6 & 12.47 & 12.31 & 12.48 & 12.64 & 97.8 &7.10 & 7.01 & 7.10 & 7.19 \\
& $\widehat{\mu}_s$ & 97.9 & 10.50 & 10.33 & 10.50 & 10.66 & 97.6 &5.59 & 5.51 & 5.59 & 5.67 \\
\hline
\end{tabular*}
\end{table}

We now examine in Table \ref{tab:confidence_bands_200}, Table \ref
{tab:confidence_bands} and Table \ref{tab:confidence_bandsc} the
empirical coverage and the width of the confidence bands, which are
built as described in Section~\ref{SCB section}.
For each sample, we estimate the covariance function $\widehat{\gamma
}_N$ and draw 10,000 realizations of a centered Gaussian process with
variance function $\widehat{\gamma}_N$ in order to obtain a suitable
coefficient $c$ with a confidence level of \mbox{$1-\alpha=0.95$} as
explained in equation (\ref{simulation threshold}). The area of the
confidence band is then
%We then compare both t empirical coverage of the different procedures
%and basic statistics on the area of the confidence bands, which is
$\int_0^T 2 c \sqrt{\widehat{\gamma}(t,t)} \, \mathrm{d}t$.
The results highlight now the interest of considering smoothing
strategies combined with the weighted cross validation bandwidth
selection criterion (\ref{def:WCV}). For stratified sampling, the use
of the unweighted cross validation criterion leads to empirical
coverage levels that are significantly below the nominal one. It also
appears that linear interpolation, which does not intend to get rid of
the noise, always gives larger confidence bands than the smoothed
estimators based on $h_{\mathrm{wcv}}$. As before, smoothing
approaches become more interesting as the number of discretization
points and the noise level increase.
The empirical coverage of the smoothed estimator is lower than the
linear interpolation estimator but remains slightly higher than the
nominal one.

% The results are very interesting and highlight the usefulness of
%smoothing when the bandwidth value is not too large. Confidence bands
%areas are always larger for linear interpolation with coverage levels
%that are always larger than the nominal levels. When considering
%strategies $h_{\textrm{mean}}/8$ and $h_{\textrm{plug}}(\widehat{
%coverage levels are still high. We can remark that $h_{
%best anymore since $h_{\textrm{mean}}/4$ has a good covering rate as
%well as smaller confidence bands. This data suggests that the optimal
%bandwidth for the mean curve estimation must be higher that of the
%covariance function, as in Kneip \textit{et al} (2001) and Cardot
%(2007).

%t10 ###
\begin{table}%[t]
\tabcolsep=4pt
\caption{(Correlated noise). Empirical coverage levels $1-\widehat
{\alpha}$ and confidence band areas for different noise levels and
bandwidth choices, with $d=400$ observation times. Units are selected
by SRSWOR or stratified sampling}
\label{tab:confidence_bandsc}
\begin{tabular*}{\textwidth}{@{\extracolsep{\fill}}llllllllllll@{}}
\hline
& & \multicolumn{5}{l}{SRSWOR} & \multicolumn{5}{l@{}}{Stratified sampling} \\[-5pt]
& & \multicolumn{5}{l}{\hrulefill} & \multicolumn{5}{l@{}}{\hrulefill} \\
$\delta$ & $h$ & $1-\widehat{\alpha}$ & Mean & 1Q & Median & 3Q &$1-\widehat{\alpha}$ & Mean & 1Q & Median & 3Q\\
\hline
5\% & lin & 97.4 & 21.33 & 21.02 & 21.32 & 21.68 & 96.9 & 5.83 & 5.75 &5.83 & 5.90 \\
& $h_\mathrm{cv}$ & 97.4 & 21.29 & 20.94 & 21.30 & 21.60 & 58.1 & 5.69& 5.61 & 5.69 & 5.77 \\
& $h_\mathrm{wcv}$ & 97.4 & 21.29 & 20.94 & 21.30 & 21.60 & 96.8 &5.79 & 5.71 & 5.79 & 5.87 \\
& $h_\mathrm{oracle}$ & 97.4 & 21.29 & 20.94 & 21.29 & 21.61 & 96.6 &5.79 & 5.71 & 5.79 & 5.87 \\
& $\widehat{\mu}_s$ & 97.4 & 20.77 & 20.42 & 20.76 & 21.10 & 97.6 &5.52 & 5.44 & 5.52 & 5.60 \\
[6pt]
25\% & lin & 98.0 & 13.51 & 13.33 & 13.52 & 13.68 & 95.7 & 7.79 & 7.71& 7.78 & 7.86 \\
& $h_\mathrm{cv}$ & 97.5 & 12.06 & 11.88 & 12.05 & 12.23 & 72.6 & 7.16& 7.05 & 7.14 & 7.24 \\
& $h_\mathrm{wcv}$ & 97.5 & 12.06 & 11.88 & 12.05 & 12.23 & 95.0 &7.53 & 7.46 & 7.53 & 7.60 \\
& $h_\mathrm{oracle}$ & 97.6 & 12.06 & 11.88 & 12.05 & 12.22 & 95.6 &7.58 & 7.50 & 7.57 & 7.66 \\
& $\widehat{\mu}_s$ & 97.2 & 10.49 & 10.31 & 10.48 & 10.66 & 97.4 &5.52 & 5.44 & 5.51 & 5.60 \\
\hline
\end{tabular*}
\end{table}

As a conclusion of this simulation study, it appears that smoothing is
not a crucial aspect when the only target is the estimation of the
mean, and that bandwidth values should be chosen by a cross validation
criterion that takes the sampling weights into account. When the goal
is also to build confidence bands, smoothing with weighted cross
validation criteria lead to narrower bands compared to interpolation
techniques, without deteriorating the empirical coverage. Smoothing
strategies which do not take account of unequal probability sampling
weights lead to empirical coverage levels that can be far below the
expected ones.

% As a conclusion of this simulation study, it appears that smoothing
%is not a crucial aspect when the only target is the estimation of the
%mean, and bandwidth values should not be chosen by usual
%cross-validation or plug-in techniques on each curve. Indeed such
%data-driven approaches select bandwidth values that are too large and
%introduce bias.
% When the goal is to build confidence bands, the situation is
%different. It appears that interpolation techniques lead to large
%widths of the bands and smoothing with intermediate bandwidth values
%can decrease in a important way the width of the band without
%deteriorating the empirical coverage.

%%%%%%%%%%%%%%%%%%%%%%%%%%%%%%%%%%%%%%%%%%%%%%%%%%%%%%%%%%%%%%%

%s5 ###
\section{Concluding remarks}\label{sec5}

% Summary of our contribution

%We have studied in this paper the use of survey sampling methods for
%estimating a population mean temporal signal.
In this paper, we have used survey sampling methods to estimate a
population mean temporal signal.
This type of approach is extremely effective when data transmission or
storage costs are important, in particular for large networks of
distributed sensors.
Considering noisy functional data, we have built
%an estimator by first smoothing the sampled curves
%and then constructing
the Horvitz--Thompson estimator of the population mean function based
on a smooth version of the sampled curves.
%based on the smoothed curves.
It has been shown that this estimator satisfies a CLT in the space of
continuous functions and that its covariance can be estimated uniformly
and consistently. Although our theoretical results were presented in
this paper with a Horvitz--Thompson covariance estimator,
they are very likely to hold for other popular estimators such as the
Sen--Yates--Grundy estimator.
We have applied our results to the construction of confidence bands
with asymptotically correct coverage.
The bands are simply obtained by simulating Gaussian processes
conditional on the estimated covariance.
%exploited our results to show that by simulating Gaussian processes
%conditional on the estimated covariance, one obtains global confidence
%bands with asymptotically correct coverage.
The problem of bandwidth selection, which is particularly difficult in
the survey sampling context, has been addressed.
We have devised a weighted cross-validation method that aims at
mimicking an oracle estimator.
This method has displayed very good performances in our numerical study;
however, a rigorous study of its theoretical properties remains to be done.
Our numerical study has also revealed that in comparison to SRSWOR,
unequal probability sampling (e.g., stratified sampling) yields far
superior performances
and that when the noise level in the data is moderate to high,
incorporating a smoothing step in the estimation procedure enhances the
accuracy in comparison to linear interpolation.
Furthermore, we have seen that even when the noise level is low,
smoothing can be beneficial for building confidence bands.
Indeed, smoothing the data leads to estimators that have higher
temporal correlation, which in turn makes the confidence bands narrower
and more stable.
Our method for confidence bands is simple and quick to implement.
It gives satisfactory coverage (a little conservative) when the
bandwidth is chosen correctly, for example, with our weighted
cross-validation method.
Such confidence bands can find a variety of applications in statistical testing.
They can be used to compare mean functions in different
sub-populations, or to test for a parametric shape or for periodicity,
among others.
Examples of applications can be found in Degras \cite{De11}.

%being able to build accurate confidence bands via procedures based on
%simulations, as proposed in this work, is an interesting and simple
%way to evaluate to which domain should belong the true unobserved mean
%trajectory.

% Comments on assumptions used in the paper and how they can be relaxed

This work also raises some questions which deserve further investigation.
A straightforward extension could be to relax the normality assumption
made on the measurement errors.
It is possible to consider more general error distributions under
additional assumptions on the moments and much longer proofs.
In another direction, it would be worthwhile to see whether our
methodology can be extended to build confidence bands for other
functional parameters such as population quantile or covariance functions.
Also, as mentioned earlier, the weighted cross-validation proposed in
this work seems a promising candidate for automatic bandwidth selection.
However, it is for now only based on heuristic arguments and its
theoretical underpinning should be investigated.

% Possible directions for future work

%As shown in our simulation study, the bandwidth is an important tuning
%parameter which has to be chosen carefully in order to get accurate
%estimations and confidence bands.

Finally, it is well known that taking account of auxiliary information,
which can be made available for all the units of the population at a
low cost,
can lead to substantial improvements with model assisted estimators (S\"
{a}rndal \textit{et al.} \cite{SaSwWr92}).
In a functional context, an interesting strategy consists in first
reducing the dimension through a functional principal components
analysis shaped for the sampling framework (Cardot \textit{et al.}
\cite{CaChGoLa10a}) and then considering semi-parametric models relating the
principal components scores to the auxiliary variables (Cardot \textit
{et al.} \cite{CaDeJo10b}). It is still possible to get consistent estimators of
the covariance function of the limit process but further investigations
are needed to prove the functional asymptotic normality and deduce that
Gaussian simulation-based approaches still lead to accurate confidence bands.

%%%%%%%%%%%%%%%%%%%%%%%%%%%%%%%%%%%%%%%%%%%%%%%%%%%%%%%%%%%%%%%
% \textbf{Acknowledgement}.
%Etienne Josserand thanks the \textit{Conseil R\'egional de Bourgogne,
%France} for its financial support (FABER PhD grant).

\begin{appendix}\label{app}

\section*{Appendix}

Throughout the proofs, we use the letter $C$ to denote a generic
constant whose value may vary from place to place.
This constant does not depend on $N$ nor on the arguments $s,t \in
[0,T]$. Note also that the expectation $\mathbb{E}$ is jointly with
respect to
the design and the multivariate normal model.

\begin{pf*}{Proof of Theorem \ref{asymptotic normality mu_hat}} %
% \vskip2mm
We first decompose the difference between the estimator $ \widehat{\mu
}_N(t)$ and its target $\mu_N(t)$ as the sum of
two stochastic components, one pertaining to the sampling variability
and the other to the measurement errors,
and of a deterministic bias component:
\renewcommand{\theequation}{\arabic{equation}}
\setcounter{equation}{21}
\begin{eqnarray}
\label{decomposition mu_hat - mu}
\widehat{\mu}_N(t) -
\mu_N(t) &=& \frac{1}{N} \sum_{k \in U}
\biggl( \frac{I_k}{\pi_k} - 1 \biggr) \tilde{X}_k(t) +
\frac{1}{N} \sum_{k \in U} \frac{I_k}{\pi_k}
\tilde{\varepsilon}_k(t) + \frac{1}{N} \sum
_k \bigl( \tilde{X}_k(t) - X_k(t)
\bigr),\qquad
\end{eqnarray}
where $\tilde{X}_k(t)$ and $ \tilde{\varepsilon}_k(t)$ are defined in
\eqref{def tilde X eps Delta}.

 \textit{Bias term}.
To study the bias term $ N^{-1} \sum_k  ( \tilde{X}_k(t) - X_k(t)
) = \mathbb{E}( \widehat{\mu}_N(t) ) - \mu_N(t) $ in \eqref
{decomposition mu_hat - mu}, it suffices to use classical results on
local linear smoothing (e.g., Tsybakov \cite{Ts09}, Proposition 1.13)
together with the H\"older continuity (A2) of the $X_k$ to see that
\begin{equation}
\label{bias bound} \sup_{t\in[0,T]} \bigg| \frac{1}{N} \sum
_k \bigl( \tilde{X}_k(t) - X_k(t)
\bigr) \bigg| \le\frac{1}{N} \sum_k
\sup_{t\in[0,T]} \bigl\llvert \tilde{X}_k(t) -
X_k(t) \bigr\rrvert \leq C h^\beta.
\end{equation}
Hence, for the bias to be negligible in the normalized estimator, it is
necessary that the bandwidth satisfy $\sqrt{N}h^{\beta} \to0$ as
$N\to\infty$.

 \textit{Error term}.
We now turn to the measurement error term %$N^{-1} \sum_{k \in U}(I_k /
in \eqref{decomposition mu_hat - mu}, which can be seen as a sequence
of random functions.
We first show that this sequence goes pointwise to zero in mean square
(a fortiori in probability) at a rate $(Ndh)^{-1}$.
We then establish its tightness in $C([0,T])$, when premultiplied by
$\sqrt{ N}$, to prove the uniformity of the convergence over $[0,T]$.

Writing the vector of local linear weights at point $t$ as
\[
W(t)=\bigl(W_1(t),\ldots,W_d(t)\bigr)'
\]
and
using the i.i.d. assumption (A4) on the $(\varepsilon_{k1},\ldots
,\varepsilon_{kd})' , k\in U_N,$
we first obtain that
\begin{eqnarray*}
\label{pointwise convergence of smooth error process} \mathbb{E} \biggl(
\frac{1}{N} \sum_{k \in U} \frac{I_k}{\pi_k} \tilde
{\varepsilon}_k(t) \biggr)^2 & = & \frac{1}{N^2} \sum
_{k \in U} \frac
{1}{\pi_k} \mathbb{E} \bigl( \tilde{
\varepsilon}_k(t) \bigr)^2
\\
&=& \frac{1}{N^2} \sum_{k \in U} \frac{1}{\pi_k}
W(t)' \mathbf {V}_N W(t).
\end{eqnarray*}
 Then, considering the facts that $\min_k
\pi_k >c $ by (A1), $\| \mathbf{V}_N\|$ is uniformly
bounded in $N$ by (A4), and exploiting a classical bound on the weights of the
local linear smoother (e.g., Tsybakov \cite{Ts09}, Lemma 1.3), we deduce
that
\begin{eqnarray}
\mathbb{E} \biggl( \frac{1}{N} \sum
_{k \in U} \frac{I_k}{\pi_k} \tilde {\varepsilon}_k(t)
\biggr)^2 &\leq& \frac{N}{(\min\pi_k)N^2}   \big\| W(t)\|^2   \big\|
\mathbf{V}_N \|
\nonumber\\[-8pt]\\[-8pt]
&\le& \frac{C}{Ndh}   .\nonumber
\end{eqnarray}

We can now prove the tightness of the sequence of processes $(N^{-1/2}
\sum_k  (I_k/\pi_k  ) \tilde{\varepsilon}_k)$. Let us
define the associated pseudo-metric\vspace*{-1pt}
\[
d^2_{\varepsilon} (s,t) = \mathbb{E} \biggl( \frac{1}{\sqrt{N}}
\sum_{k \in U} \frac{I_k}{\pi_k} \bigl( \tilde{
\varepsilon}_k(s)- \tilde {\varepsilon}_k(t) \bigr)
\biggr)^2 .
\]
We use the following maximal inequality holding for sub-Gaussian processes
(van der Vaart and Wellner \cite{VaWe00}, Corollary 2.2.8):
\begin{equation}
\label{maximal inequality subgaussian process} \mathbb{E} \biggl(
\sup_{t\in[0,T]} \biggl\llvert \frac{1}{\sqrt{N}} \sum
_{k \in
U} \frac{I_k}{\pi_k} \tilde{\varepsilon}_k(t)
\biggr\rrvert \biggr) \le \mathbb{E} \biggl( \biggl\llvert \frac{1}{\sqrt{N}} \sum
_{k \in U} \frac
{I_k}{\pi_k} \tilde{\varepsilon}_k(t_0)
\biggr\rrvert \biggr) + K \int_0^\infty \sqrt{
\log N(x,d_\varepsilon)} \,\mathrm{d}x   ,
\end{equation}
where $t_0$ is an arbitrary point in $[0,T]$ and the covering number
$N(x,d_\varepsilon)$ is the minimal number of $d_\varepsilon$-balls of
radius $x>0$ needed to cover $[0,T]$. Note the equivalence of working
with packing or covering numbers in maximal inequalities, see \textit{ibid} page 98.
%(The use of covering numbers in the previous maximal inequality is
%equivalent to that of packing numbers, see van der Vaart (2000) p. 98.)
Also note that the sub-Gaussian nature of the smoothed error process
$N^{-1/2} \sum_{k \in U}  (I_k/\pi_k ) \tilde{\varepsilon}_k$
stems from the i.i.d. multivariate normality of the random vectors
$(\varepsilon_{k1},\ldots, \varepsilon_{kd})' $ and the boundedness of the $I_k$
for $k\in U_N$.

By the arguments used in \eqref{pointwise convergence of smooth error process}
and an elementary bound on the increments of the weight function vector $W$
%relying on the Lipschitz continuity of the kernel $K$
(see, e.g., Lemma 1 in Degras \cite{De11}),
one obtains that
\begin{eqnarray}
\label{entropy calculation 1} d^2_{\varepsilon} (s,t) & =&
\frac{1}{N} \sum_{k \in U} \frac{1}{\pi_k}
\mathbb{E} \bigl( \tilde{\varepsilon}_k(s)- \tilde{\varepsilon
}_k(t) \bigr)^2\nonumber
\\
& \le&\frac{1}{\min\pi_k} \big\| W(s) - W(t) \big\|^2  \|
\mathbf{V}_N \|
\\
& \le&\frac{C}{dh} \biggl( \frac{|s-t|^2}{h^2} \wedge1 \biggr) .\nonumber
\end{eqnarray}

It follows that the covering numbers satisfy
\[
\cases{
 N(x,d_\varepsilon) = 1, & \quad  if $\dfrac{C}{dh} \le x^2$,
\cr
N(x,d_\varepsilon) \le\dfrac{\sqrt{C}}{ h \sqrt{dh} x } , &\quad   if $
\dfrac{C}{dh} > x^2$. } %
\]

Plugging this bound and the pointwise convergence \eqref{pointwise
convergence of smooth error process} in the maximal inequality \eqref
{maximal inequality subgaussian process}, we get after a simple
integral calculation (see equation (17) in Degras \cite{De11} for details) that
\begin{equation}
\mathbb{E} \biggl( \sup_{t\in[0,T]} \biggl\llvert \frac{1}{\sqrt{N}} \sum
_{k \in
U} \frac{I_k}{\pi_k} \tilde{
\varepsilon}_k(t) \biggr\rrvert \biggr) \le\frac{C}{dh} + C \sqrt{
\frac{|\log(h)|}{dh}}   .
\end{equation}
Thanks to Markov's inequality, the previous bound guarantees the
uniform convergence in probability of $N^{-1/2} \sum_{k \in U}
(I_k/\pi_k) \tilde{\varepsilon}_k$ to zero, provided that $|\log
(h)|/(dh)\to0$ as $N\to\infty$. The last condition is equivalent to
$\log(d)/(dh)\to0$ by the fact that $dh\to\infty$ and
by the properties of the logarithm.
%$\lim_{x\to\infty} (\log x)/x =0 $.

 \textit{Main term: sampling variability}.
Finally, we look at the process $ N^{-1} \sum_{k \in U}  (
I_k/\pi_k - 1 ) \tilde{X}_k$
in \eqref{decomposition mu_hat - mu}, which is asymptotically normal
in $C([0,T])$ as we shall see.
We first establish the finite-dimensional asymptotic normality of this
process normalized by $\sqrt{N}$,
after which we will prove its tightness thanks to a maximal inequality.

Let us start by verifying that the limit covariance function of the
process is indeed the function $\gamma$ defined in Section \ref{AN HT
estimator}.
The finite-sample covariance function is expressed
\begin{eqnarray}
\label{limit cov mu_hat}
&&\mathbb{E} \biggl\{ \biggl(
\frac{1}{\sqrt{N}} \sum_{k \in U} \biggl(
\frac
{I_k}{\pi_k} - 1 \biggr) \tilde{X}_k(s) \biggr) \biggl( \frac{1}{\sqrt{N}} \sum_{l \in U}
\biggl( \frac{I_l}{\pi_l} - 1 \biggr) \tilde{X}_l(t) \biggr) \biggr
\}
\nonumber\\
&&\quad  = \frac{1}{N} \sum_{k,l \in U} \frac{\Delta_{kl}}{\pi_k\pi_l}
\tilde{X}_k(s)\tilde{X}_l(t)
\\
&&\quad  = \frac{1}{N} \sum_{k,l \in U} \frac{\Delta_{kl}}{\pi_k\pi_l}
{X}_k(s) {X}_l(t) + \mathcal{O} \bigl(h^\beta
\bigr)
\nonumber
\\
&&\quad  = \gamma(s,t) +o(1) + \mathcal{O} \bigl( h^\beta \bigr) .\nonumber
\end{eqnarray}
To derive the previous relation, we have used the facts that
\[
\max_{k,l\in U}\sup_{
s,t\in[0,T]} \bigl\llvert \tilde{X}_k(s)
\tilde{X}_l(t) - X_k(s)X_l(t) \bigr\rrvert
\le C h^\beta
\]
%
%since $\sup_{k,t} | \tilde{X}_k (t) - X_k(t) | =\mathcal{O} (h^\beta
%)$ and $\sup_{k,t} | \tilde{X}_k (t) | =\mathcal{O}(1) $
by \eqref{bias bound} and the uniform boundedness of the $X_k$ arising
from (A2) and that, by (A1),
\begin{eqnarray}
\label{bound delta_kl} \frac{1}{N} \sum
_{k , l\in U} \frac{|\Delta_{kl}|}{\pi_k\pi_l} & = &\frac{1}{N} \sum
_{k \ne l} \frac{|\Delta_{kl}|}{\pi_k\pi_l} + \frac{1}{N} \sum
_{k} \frac{\Delta_{kk}}{\pi_k^2}
\nonumber\\[-8pt]\\[-8pt]
& \le&\frac{1}{N} \frac{N(N-1)}{2} \frac{\max_{k,l}(n|\Delta_{kl}|)}{n} + \frac{1}{N}
\sum_{k} \frac{1-\pi_k}{\pi_k} \le C .\nonumber
\end{eqnarray}

We now check the finite-dimensional convergence of $N^{-1/2} \sum_{k
\in U}  ( I_k/\pi_k - 1 ) \tilde{X}_k$
to a centered Gaussian process with covariance $\gamma$.
In light of the Cramer--Wold theorem, this convergence is easily shown
with characteristic functions and appears as a straightforward
consequence of (A5). It suffices for us to check that the uniform
boundedness of the trajectories $X_k$ derived from (A2)
is preserved by local linear smoothing, so that the $\tilde{X}_k$ are
uniformly bounded as well.

It remains to establish the tightness of the previous sequence of
processes so as to obtain its asymptotic normality in $C([0,T])$.
To that intent we use the maximal inequality of the Corollary~2.2.5 in
van der Vaart and Wellner \cite{VaWe00}.
With the notations of this result, we consider the pseudo-metric
$d^2_{\tilde{X}}(s,t)=\mathbb{E}\{ N^{-1/2} \sum_{k \in U} ( I_k/\pi_k - 1
) ( \tilde{X}_k(s) -\tilde{X}_k(t) ) \}^2 $ and the function $\psi
(t)=t^2$ for the Orlicz norm. We get the following bound for the
second moment of the maximal increment:
\begin{eqnarray}
\label{max inequality 2} %
&&\mathbb{E} \biggl\{
\sup_{d_{\tilde{X}}(s,t)\le\delta} \biggl\llvert \frac
{1}{\sqrt{N}} \sum
_{k \in U}  \biggl(
\frac{ I_k}{\pi_k} - 1 \biggr) \bigl( \tilde{X}_k(s) - \tilde
{X}_k(t) \bigr) \biggr\rrvert \biggr\}^2
\nonumber\\[-8pt]\\[-8pt]
&&\quad  \le C \biggl( \int_0^{\eta} \psi^{-1}
\bigl( N(x,d_{\tilde{X}})\bigr) \,\mathrm{d}x + \delta\psi^{-1}
\bigl(N^2(\eta,d_{\tilde{X}})\bigr) \biggr)^2\nonumber
\end{eqnarray}
for any arbitrary constants $\eta,\delta>0$.
% $D_T= \sup_{s,t\in[0,T]}d(s,t)$ is the diameter of $[0,T]$.
Observe that the maximal inequality \eqref{max inequality 2} is weaker than
\eqref{maximal inequality subgaussian process} where an additional
assumption of sub-Gaussianity is made (no $\log$ factor in the
integral above).
Employing again the arguments of \eqref{limit cov mu_hat}, we see that
\begin{eqnarray}
d^2_{\tilde{X}}(s,t) & =& \frac{1}{N} \sum
_{k,l} \frac{\Delta_{kl}}{\pi_k\pi_l} \bigl( \tilde{X}_k(s)-
\tilde{X}_k(t) \bigr) \bigl( \tilde{X}_l(s)-
\tilde{X}_l(t) \bigr)
\nonumber
\\
& \le&\frac{C}{N}   \frac{N(N-1)}{2n} |s-t|^{2\beta} +
\frac
{C}{N}  N |s-t|^{2\beta}
\\
& \le& C |s-t|^{2\beta} .\nonumber
\end{eqnarray}
It follows that the covering number satisfies $N(x,d_{\tilde{X}})\le C
x^{-1/\beta}$
and that the integral in \eqref{max inequality 2} is smaller than $C
\int_0^\eta x^{-0.5/\beta} \,\mathrm{d}x = C\eta^{1-0.5/\beta} $,
which can be made arbitrarily small since $\beta>0.5$. Once $\eta$ is
fixed, $\delta$ can be adjusted to make the other term
in the right-hand side of \eqref{max inequality 2} arbitrarily small
as well.
With Markov's inequality, we deduce that the sequence $(N^{-1/2} \sum_{k \in U}  ( I_k/\pi_k - 1 ) \tilde{X}_k)_{N\ge1}$
is asymptotically $d_{\tilde{X}}$-equicontinuous in probability (with
the terminology of van der Vaart and Wellner \cite{VaWe00}), which guarantees
its tightness in $C([0,T])$.
\end{pf*}

%%%%%%%%%%%%%%%%%%%%%%%%%%

\begin{pf*}{Proof of Theorem \ref{teo2}}
To establish the uniform convergence of the covariance estimator, we
first show its mean square convergence in the pointwise sense.
Then, we extend the pointwise convergence to uniform convergence
through an asymptotic tightness argument (i.e., by showing that for
$N$ large enough, the covariance estimator lies in a compact $K$ of
$C([0,T]^2)$ equipped with the sup-norm with probability close to 1).
We make use of maximal inequalities to prove the asymptotic tightness result.
%of the covariance estimator

 \textit{Mean square convergence}.
%We start by showing that
We first decompose the distance between $ \widehat{\gamma}_N(s,t)$
and its target $ \gamma_N(s,t)$
as follows:
%goes to zero in mean square as $N\to\infty$.
%
\renewcommand{\theequation}{\arabic{equation}}
\setcounter{equation}{31}
\begin{eqnarray}
\label{decomposition gamma_hat} \widehat{\gamma}_N(s,t)
- \gamma_N(s,t) & =& \frac{1}{N} \sum
_{k,l\in U} \frac{\Delta_{kl}}{\pi_k\pi_l} \biggl( \frac{I_k I_l}{\pi_{kl}} -1 \biggr)
\tilde{X}_{k}(s) \tilde{X}_{l} (t)
\nonumber
\\
&&{} + \frac{1}{N} \sum_{k,l\in U}
\frac{\Delta_{kl}}{\pi_k\pi_l} \frac{I_k I_l}{\pi_{kl}} \bigl( \tilde{X}_{k}(s) \tilde {
\varepsilon}_{l} (t) + \tilde{X}_{l}(t) \tilde{
\varepsilon}_{k} (s) \bigr)
\nonumber
\\
&&{} + \frac{1}{N} \sum_{k,l\in U}
\frac{\Delta_{kl}}{\pi_k\pi_l} \frac{I_k I_l}{\pi_{kl}} \tilde{\varepsilon}_{k}(s) \tilde {
\varepsilon}_{l} (t)
\\
& &{} - \frac{1}{N} \sum_{k \in U}
\frac{1}{\pi_k} \mathbb {E} \bigl( \tilde{\varepsilon}_{k}(s) \tilde{
\varepsilon}_{k} (t) \bigr)
\nonumber
\\
& :=& A_{1,N} + A_{2,N} + A_{3,N}
-A_{4,N}.\nonumber
\end{eqnarray}

To establish the mean square convergence of $ ( \widehat{\gamma
}_N(s,t)- \gamma_N(s,t) ) $ to zero as
$N\to\infty$, it is enough to show that $\mathbb{E}( A_{i,N}^2) \to
0$ for
$ i=1,\ldots,4,$ by the Cauchy--Schwarz inequality.

%&\mathbb{E}\big( \widehat{\gamma}_N(s,t) - \gamma_N(s,t)\big)^2 \\
%& = \frac{1}{N^2} \sum_{k,l}\sum_{k',l'} \left( \frac{\Delta_{kl}
%I_{l'}}{\pi_{k'l'}} -1 \right)\widetilde{X}_{k}(s) \widetilde{X} _{l}
%&+ \frac{1}{N^2} \sum_{k,l}\sum_{k',l'} \left( \frac{\Delta_{kl}
%& := A_N + B_N.
Let us start with
\begin{eqnarray}
\label{A1n}
\mathbb{E}\bigl(A_{1,N}^2\bigr)
&=& \frac
{1}{N^2} \sum_{k,l}\sum
_{k',l'} \frac{\Delta_{kl} \Delta_{k'l'} }{ \pi_k\pi_l
\pi_{k'}\pi_{l'} }   \frac{ \mathbb{E} \{ (I_{k}I_{l} - \pi_{kl})(I_{k'}I_{l'} - \pi_{k'l'}) \}}{\pi_{kl}\pi_{k'l'}}
\nonumber\\[-8pt]\\[-8pt]
&&\hspace*{47pt}{}\times\widetilde{X}_{k}(s) \widetilde{X}_{l}(t)
\widetilde{X}_{k'}(s) \widetilde{X}_{l'}(t) .\nonumber
\end{eqnarray}
It can be shown that this sum converges to zero by strictly following
the proof of the Theorem~3 in Breidt and Opsomer \cite{BrOp00}. The idea of
the proof is to partition the set of indexes in \eqref{A1n} into (i)~$k=l$ and $k'=l'$, (ii) $k=l$ and $k'\ne l'$ or vice-versa, (iii) $k\ne
l$ and $k' \ne l'$, and study the related subsums. The convergence to zero
is then handled with assumption (A1) (mostly) in case (i), with
(A1)--(A6) in case (iii), and thanks to the previous results and
Cauchy--Schwarz inequality in case (ii). More precisely, it holds that
\begin{eqnarray}
  \mathbb{E}\bigl(A_{1,N}^2
\bigr) &\le&\frac{C\max_{k\ne l} n|\Delta_{kl}|}{(\min\pi_k)^4 n} + \frac{C}{(\min\pi_k)^3N}
\nonumber\\
& &{} + \biggl( \frac{C(\max_{k\ne l} n|\Delta_{kl}|)N}{(\min\pi_k)^2 (\min_{k\ne l} \pi_{kl}) n } \biggr)^2
\\
&&\hspace*{12pt}{}\times
\max_{(k,l,k',l') \in D_{4,N}} \big| \mathbb{E} \bigl\{ (I_{k}I_{l} -
\pi_{kl}) (I_{k'}I_{l'} - \pi_{k'l'})
\bigr\} \big| .\qquad \nonumber
\end{eqnarray}

For the (slightly simpler) study of $ \mathbb{E}(A_{2,N}^2)$, we
provide an
explicit decomposition:
\begin{eqnarray}
 \mathbb{E}\bigl(A_{2,N}^2\bigr)
& = &\frac{4}{N^2} \sum_{k, l} \sum
_{k' } \frac
{\Delta_{kl}\Delta_{k'l} }{ \pi_k \pi_{k'}\pi_l^2 } \tilde {X}_{k}(s)
\tilde{X}_{k'}(t)   \mathbb{E}( I_kI_{k'}I_l)
  \mathbb{E} \bigl( \tilde{\varepsilon }_{l}(s) \tilde {
\varepsilon}_{l}(t) \bigr)\nonumber
\\
& = &\frac{4}{N^2} \sum_{k\in U}
\frac{\Delta_{kk}^2 }{ \pi_k^5 } \tilde{X}_{k}(s) \tilde{X}_{k}(t)
\mathbb{E} \bigl( \tilde{\varepsilon}_{k}(s) \tilde{\varepsilon}_{k}(t)
\bigr)
\nonumber\\[-8pt]\\[-8pt]
%&  + \frac{4}{N^2} \sum_{k\ne k'} \frac{\Delta_{kk}\Delta_{kk'}
%}{ \pi_k^4 \pi_{k'} } \tilde{X}_{k}(s) \tilde{X} _{k'}(t)  \E(
& &{} + \frac{8}{N^2} \sum
_{k\ne k'} \frac{\Delta_{kk}\Delta_{kk'} }{ \pi_k^4 \pi_{k'} \pi_{kk'}} \tilde{X}_{k}(s)
\tilde{X}_{k'}(t)  \mathbb{E}\bigl( \tilde{\varepsilon}_{k}(s)
\tilde{\varepsilon}_{k}(t) \bigr)\nonumber
\\
% &  + \frac{4}{N^2} \sum_{k\ne l} \sum_{k': k' \ne l} \frac{
& &{} + \frac{4}{N^2} \sum
_{k,k'} \sum_{ l \notin\{k,k'\}}
\frac
{\Delta_{kl}\Delta_{k'l} }{ \pi_k \pi_{k'}\pi_l^2 \pi_{kl}\pi_{k'l}}   \tilde{X}_{k}(s) \tilde{X}_{k'}(t)
\mathbb{E}( I_kI_{k'}I_l)   \mathbb{E} \bigl(
\tilde{\varepsilon}_{l}(s) \tilde{\varepsilon}_{l}(t) \bigr) .\nonumber
\end{eqnarray}
Note that the expression of $ \mathbb{E}(A_{2,N}^2) $ as a quadruple
sum over
$k,l,k',l' \in U_N$ reduces to a triple sum since $\mathbb{E}( \tilde
{\varepsilon}_{l}(s) \tilde{\varepsilon}_{l'}(t) ) = 0$ if $l\ne l'$ by
%%independence of the error vectors $(\varepsilon_{k1},\ldots,
(A4).
Also note that $| \mathbb{E}( I_kI_{k'}I_l) | \le1$ for all
$k,k',l\in U$.
With (A1), (A2), and the bound $ |\mathbb{E}( \tilde{\varepsilon}_{k}(s)
\tilde{\varepsilon}_{k}(t) ) |=| W(s)'\mathbf{V}_N W(t)| \le\| W(s)
\| \| \mathbf{V}_N \| \| W(t)\| \le C/(dh)$,
it follows that
\begin{eqnarray}
\mathbb{E}\bigl(A_{2,N}^2\bigr) &\le&\frac{CN}{N^2}
\frac{\| \mathbf{V}_N \|
}{dh} + \frac{CN^2}{N^2} \frac{\max_{k\ne k'} n|\Delta_{kk'}| }{n} \frac{\| \mathbf{V}_N \|
}{dh}
\nonumber\\[-8pt]\\[-8pt]
& &{} + \frac{CN^3}{N^2} \frac{(\max_{k\ne l} n|\Delta_{kl}| )^2
}{n^2} \frac{\| \mathbf{V}_N \|}{dh} =
\frac{C}{Ndh}.\nonumber
\end{eqnarray}

% To study the term $ \E(A_{3,N}^2)$, we start with the same partition
%of the quadruple sum as the one used with $ \E(A_{1,N}^2)$.
%Here, due to the independence (A4) of the error vectors,
% the partition simplifies further into (i) $k=l$, $k'=l'$, $k\ne k'$,
%and (ii) $k=l=k'=l'$:
We turn to the evaluation of
\[
\mathbb{E}\bigl(A_{3,N}^2\bigr)= \frac{1}{N^2} \sum
_{k,l,k',l'} \frac{\Delta_{kl}\Delta_{k'l'}}{\pi_k\pi_l\pi_{k'}\pi_{l'}} \frac{\mathbb{E}(I_k
I_lI_{k'}I_{l'})}{\pi_{kl}\pi_{k'l'}}
\mathbb{E} \bigl( \tilde {\varepsilon }_{k}(s) \tilde{
\varepsilon}_{l} (t) \tilde{\varepsilon}_{k'}(s) \tilde {
\varepsilon}_{l'} (t) \bigr).
\]
We use the independence (A4) of the errors across population units
to partition
the above quadruple sum $ \mathbb{E}(A_{3,N}^2)$ according to the cases
(i) $k=l$, $k'=l'$, $k\ne k'$, (ii) $k=l'$, $k'=l$, and $k\ne k'$,
(iii) $k=k'$, $l=l'$, and $k\ne l$
and (iv) $k=l=k'=l'$. Therefore,
\begin{eqnarray}
\label{variance A3n} \mathbb{E}\bigl(A_{3,N}^2\bigr)& =&
\frac{1}{N^2} \sum_{k \ne k' } \frac{\pi_{kk'}}{\pi_k^2\pi_{k'}^2}
\biggl( \frac{\Delta_{kk}\Delta_{k'k'}}{\pi_k\pi_{k'}} + \frac
{\Delta_{kk'}^2}{\pi_{kk'}^2} \biggr) \mathbb{E}\bigl(\tilde{
\varepsilon}_{k}(s) \tilde{\varepsilon}_{k} (t)\bigr)
\mathbb{E} \bigl(\tilde{\varepsilon}_{k'}(s) \tilde{\varepsilon}_{k'}
(t) \bigr)
\nonumber\\[-8pt]\\[-8pt]
&&{} + \frac{1}{N^2} \sum_{k \ne l }
\frac{\Delta_{kl}^2}{\pi_k^2\pi_l^2\pi_{kl}} \mathbb{E}\bigl( \tilde{\varepsilon}_{k}^2
(s)\bigr)   \mathbb{E}\bigl(\tilde {\varepsilon}_{l}^2(t)
\bigr) % &  + \frac{1}{N^2} \sum_{ k } \frac{\Delta_{kk}}{\pi_k^2}
%_{k}^2 (t) ) .
+ \frac{1}{N^2} \sum
_{ k } \frac{\Delta_{kk}^2}{\pi_k^5}  \mathbb{E} \bigl(\tilde{
\varepsilon}_{k}^2(s) \tilde{\varepsilon}_{k}^2
(t) \bigr) .\nonumber
\end{eqnarray}

Forgoing the calculations already done before, we focus on the main
task which for this term is to bound the quantity $\mathbb{E}(\tilde
{\varepsilon
}_{k}^2(s)\tilde{\varepsilon}_{k}^2(t) )$ (recall that $\mathbb
{E}(\tilde
{\varepsilon}_{k}(s)\tilde{\varepsilon}_{k}(t)) \le C / (dh)$ as seen before).
We first note that $\mathbb{E}(\tilde{\varepsilon}_{k}^2(s)\tilde
{\varepsilon
}_{k}^2(t)) \le \{\mathbb{E}(\tilde{\varepsilon}_{k}^4(s))
\}^{1/2}  \{\mathbb{E}(\tilde{\varepsilon}_{k}^4(t))  \}^{1/2}.$ %
Writing $\bolds{\varepsilon} \sim N(0,\mathbf{V}_N)$,
it holds that $ \mathbb{E}(\tilde{\varepsilon}_{k}^4(t)) = \mathbb{E}(
(W(t)'\bolds{\varepsilon})^4 ) = 3 (W(t)' \mathbf{V}_N W(t) )^2$
by the moment properties of the normal distribution. Plugging this
expression in \eqref{variance A3n},
we find that
\begin{equation}
\mathbb{E}\bigl(A_{3,N}^2\bigr) \le\frac{C}{(dh)^2} +
\frac{C}{ N(dh)^{2}} .
\end{equation}

Finally, like $\mathbb{E}(\tilde{\varepsilon}_{k}(s)\tilde{\varepsilon}_{k}(t))$,
the deterministic term $A_{4,N}$ is of order $1/(dh)$.
%by following the argument of \eqref{pointwise convergence of smooth
%error process}.
%One easily concludes that. % uniformly in $s,t\in[0,T]$.

%w is problem has been studied before
%The quadruple sum $B_N$ can be studied by
%decomposing the summands according to the cardinality of
%$ \big\{k,l , k',l' \big\}$. % which belongs to $ \{ 0,2,4\}$.
%When this cardinality is 3 or 4, meaning that at least one index is
%different from all three others, the expectation in the summand is
%null by virtue of the independence between the sampling and the errors
%and of the independence among the errors $(\varepsilon_{k1},\ldots,
%e.g. when $k=l$ and $k'=l'$ but $k\ne k'$, the expectation in the
%summand reduces to
%a product of covariances of error terms, again by (A5) ($\operatorname{Cov}
%_{k'}(t)\big) $ in our example). These covariances have already been
%studied in \eqref{entropy calculation 1} and are of order $(ph)^{-1}$.
%Now, there are ${4\choose2} {N \choose2} \asymp N^2$ subsets $
%of cardinality 2 and in this case the ratios $ \Delta_{kl}
%This guarantees
%that the sum over these subsets goes indeed to zero, when
%premultiplied by $N^{-2}$.
%It remains to study the case where $ \# \big\{k,l , k',l' \big\} =1$.
%In this case, the expectation in the summand reads $\operatorname{Var} \big(
%us write $\operatorname{Var} \big( \widetilde{\varepsilon} _{k}(s) \widetilde{
%_{k}(s) \widetilde{\varepsilon} _{k}(t)\big)\big]^2 $.
%The second term $ \big[ \E( \widetilde{\varepsilon} _{k}(s) \widetilde{

 \textit{Tightness}.
To prove the tightness of the sequence $ (\widehat{\gamma}_N- \gamma_N )_{ N\ge1}$ in $C([0,T]^2)$,
we study separately each term in the decomposition \eqref
{decomposition gamma_hat} and we
call again to the maximal inequalities of van der Vaart and Wellner \cite{VaWe00}.

For the first term $A_{1,N}=A_{1,N}(s,t)$, we consider the
pseudo-metric $d$ defined as the $L^4$-norm of the increments:
$d^4_1((s,t),(s',t')) = \mathbb{E}| A_{1,N}(s,t) - A_{1,N}(s',t') |^4
$. (The
need to use here the $L^4$-norm and not the usual $L^2$-norm
is justified hereafter by a dimension argument.)
With \mbox{(A1)--(A2)} and the approximation properties of local linear
smoothers, one sees that
\[
\biggl\llvert \frac{1}{N} \sum_{k,l\in U}
\frac{\Delta_{kl}}{\pi_k\pi_l} \biggl( \frac{I_k I_l}{\pi_{kl}} -1 \biggr) \bigl( \tilde
{X}_{k}(s) \tilde{X}_{l} (t) - \tilde{X}_{k}
\bigl(s'\bigr) \tilde{X}_{l} \bigl(t'\bigr)
\bigr) \biggr\rrvert \le C \bigl( \big|s-s'\big|^\beta+
\big|t-t'\big|^\beta \bigr).
\]
Hence $d_1(s,t)$ $\le$ $ C  ( |s-s'|^\beta+ |t-t'|^\beta
) $ and for all $x>0$,
the covering number $N(x,d_1)$ is no larger than the size of a
two-dimensional square grid of mesh $x^{1/\beta}$, that is,
$N(x,d_1)$ $ \le$ $ C x^{-2/\beta}$.
(Compare to the proof of Theorem \ref{asymptotic normality mu_hat} where, for the main term $N^{-1/2}
\sum_k (I_k/\pi_k) \tilde{X}_k$, we have
$N(x,d_{\tilde{X}})\le C x^{-1/\beta}$ because the index set $[0,T]$
is of dimension 1.)
Using Theorem 2.2.4 of van der Vaart and Wellner \cite{VaWe00} with $\psi
(t)=t^4$, it follows that for all $\eta,\delta>0,$
\begin{eqnarray*}
&&\mathbb{E} \Bigl\{ \sup_{d_1((s,t),(s',t'))\le\delta} \big| A_{1,N}(s,t) -
A_{1,N}\bigl(s',t'\bigr)\big|^4
\Bigr\}\\
&&\quad  \le C \biggl( \int_{0}^\eta
\psi^{-1}\bigl(N(x,d_1)\bigr)\,\mathrm{d}x + \delta\psi^{-1}
\bigl(N^2(\eta,d_1)\bigr) \biggr)^4
\\
&&\quad  \le C \bigl( \eta^{1-0.5/\beta} + \delta\eta^{-1/\beta}
\bigr)^4 .
\end{eqnarray*}
The upper bound above can be made arbitrarily small by varying $\eta$
first and $\delta$ next since $\beta> 0.5$. Hence, with
Markov's inequality, we deduce that the processes $A_{1,N}$
are tight in $C([0,T]^2)$. %(Note that these processes are sub-Gaussian
%as they are weighted -functional- averages of bounded random
%variables.)

The bivariate processes $ (A_{2,N})_{N\ge1} $ are sub-Gaussian for the
same reasons as the univariate processes
$N^{-1/2} \sum_{k \in U}  (I_k/\pi_k ) \tilde{\varepsilon
}_k$ are in the proof of Theorem \ref{asymptotic normality mu_hat}, namely
the independence and multivariate normality of the error vectors
$(\varepsilon_{k1},\ldots, \varepsilon_{kd})' $
and the boundedness of the sample membership indicators $I_k$ for $k\in
U_N$. Therefore,
although the covering number $N(x,d_2)$ grows to $\mathrm{O}(x^{-2/\beta})$ in
dimension 2,
with $d_2$ being the $L^2$-norm on $[0,T]^2$, this does not affect
significantly the integral upper bound
$\int_{0}^\infty\sqrt{\log(N(x,d_2))}\,\mathrm{d}x $ in a maximal inequality
like \eqref{maximal inequality subgaussian process}.
As a consequence, one obtains the tightness of $(A_{2,N})$ in $C([0,T]^2)$.

To study the term $A_{3,N}(s,t)$ in \eqref{decomposition gamma_hat},
we start with the following bound:
% transform
%
\begin{eqnarray*}
\big| A_{3,N}(s,t) \big| & \le&\frac{1}{N} \sum
_{k,l} \frac
{|\Delta_{kl}|}{\pi_k\pi_l} \frac{I_kI_l}{\pi_{kl}}
\frac{\tilde
{\varepsilon}^2_k(s)+ \tilde{\varepsilon}_l^2(t)}{2}
\\
& =& \frac{1}{N} \sum_k \biggl( \sum
_l \frac{|\Delta_{kl}| }{2\pi_l}\frac{I_l}{\pi_{kl}} \biggr)
\frac{I_k}{\pi_k} \tilde{\varepsilon }^2_k(s) +
\frac{1}{N} \sum_l \biggl( \sum
_k \frac{|\Delta_{kl}| }{2\pi_k}\frac{I_k}{\pi_{kl}} \biggr)
\frac{I_l}{\pi_l} \tilde{\varepsilon }^2_l(t)
\\
& \le&\frac{C}{N} \sum_{k} \tilde{
\varepsilon}^2_k(s) + \frac{C}{N} \sum
_{l} \tilde{\varepsilon}^2_l(t) .
%& = \frac{C}{N} \left( W(s)' W_N W(s) + W(t)' W_N W(t) \right)
\end{eqnarray*}

%where the $p\times p$ random matrix $ W_N$ has a Wishart distribution
%with parameters $\mathbf{V}_N$ and $N$ (degrees of freedom).
%It suffices now to check that the process $(N^{-1} W(t)' W_N W(t))_{t
%Using the a classical limit result for the largest eigenvalue of a
%Wishart matrix (e.g. Johnstone (2001))
%of

  The two-dimensional study is thus reduced to an easier
one-dimensional problem.

To apply the Corollary 2.2.5 of van der Vaart and Wellner \cite{VaWe00}, we
consider the function $\psi(t)=t^m$ and
the pseudo-metric $d_3^m(s,t)= \mathbb{E}| N^{-1} \sum_{k} ( \tilde
{\varepsilon
}^2_k(s)- \tilde{\varepsilon}^2_k(t) ) |^m$,
where $m\ge1$ is an arbitrary integer.
We have that
\begin{equation}
\label{coro 2.2.5 VW2000} \mathbb{E} \biggl\{ \sup_{s,t\in[0,T]} \biggl\llvert
\frac{1}{N} \sum_{k} \bigl( \tilde{
\varepsilon}^2_k(s)-\tilde{\varepsilon}^2_k(t)
\bigr)\biggr\rrvert^m \biggr\} \le C \biggl( \int
_{0}^{D_T} \bigl(N(x,d_3)
\bigr)^{1/m}\,\mathrm{d}x \biggr)^m,
\end{equation}
where $D_T = \sup_{s,t\in[0,T]}d_3(s,t)$ is the diameter of $[0,T]$
for $d_3$.
Using the classical inequality, $\llvert  \sum_{k=1}^n a_k \rrvert^m
\leq n^{m-1} \sum_{k=1}^n \llvert  a_k \rrvert^m,$ for $m>1$ and
arbitrary real numbers $a_1, \ldots, a_n,$ we get, with the
Cauchy--Schwarz inequality and the moment properties of Gaussian random
vectors, that
\begin{eqnarray}
\label{d for chi square process} d_3^m(s,t) & \leq&
\frac{1}{N} \sum_k \mathbb{E}\bigl\llvert
\tilde {\varepsilon }^2_k(s)-\tilde{\varepsilon}^2_k(t)
\bigr\rrvert^m
\nonumber
\\
& \leq&\frac{1}{N} \sum_k \bigl\{\mathbb{E}
\bigl\llvert \tilde{\varepsilon }_k(s)-\tilde{\varepsilon}_k(t)
\bigr\rrvert^{2m} \bigr\}^{1/2} \bigl\{ \mathbb{E} \bigl
\llvert \tilde{\varepsilon}_k(s)+\tilde{\varepsilon}_k(t)
\bigr\rrvert^{2m} \bigr\}^{1/2}
\\
& \leq&\frac{C_m}{N} \sum_k \bigl\llVert
W(s)-W(t)\bigr\rrVert_{\mathbf
{V}_N}^{m} \bigl\llVert W(s)+W(t)\bigr
\rrVert_{\mathbf{V}_N}^{m}
\leq\frac{C_m'}{(dh)^m} \biggl( \frac{|s-t|}{h} \wedge1
\biggr)^m,\nonumber
\end{eqnarray}
where $\| \mathbf{x} \|_{\mathbf{V}_N} = (\mathbf{x'V}_N\mathbf
{x})^{1/2}$ and $C_m$ and $C_m'$ are constants that only depend on $m.$

%Let $\mathcal{I}_{m,N} = \{ i =(i_1,\ldots,i_m): 1\le i_1 \le\cdots
%be the set of all (ordered but not necessarily with distinct elements)
%$m$-uples of $U_N$.
%For each $i\in\mathcal{I}_{m,N}$, let $r(i)$ be the number of
%distinct elements and $G_{i1},\ldots, G_{ir(i)}$ the associated
%partition of $i$
%in groups of indexes pointing to the same element in $U_N$.
%and $a \wedge b = \min(a,b)$. %We observe that the resulting bound
%for $d(s,t)$ is quite similar to \eqref{entropy calculation 1}.

We deduce from \eqref{d for chi square process} that the diameter
$D_T$ is at most of order $1/(dh)$
and that for all $ 0<x \le1/(dh)$, the covering number $N(x,d_3)$ is
of order $1/(xdh^2)$.
Hence, the integral bound in \eqref{coro 2.2.5 VW2000}
is of order $\int_{0}^{1/(dh)} (dh^2x)^{-1/m} \,\mathrm{d}x \le C (dh^2)^{-1/m}
(dh)^{(1-1/m)} = C/(dh)^{1+1/m}$.
Therefore, if $dh^{1+\alpha} \to\infty$ for some $\alpha>0,$ the sequence
$( N^{-1} \sum_{k} ( \tilde{\varepsilon}^2_k))_{N\ge1}$ tends
uniformly to zero in probability which concludes the study of the term
$(A_{3,N})_{N\ge1}$ and the proof.
\end{pf*}

\begin{pf*}{Proof of Theorem \ref{conditional weak functional
convergence for Gaussian processes}}
We show here the weak convergence of $(\widehat{G}_N)$ to $G$ in
$C([0,T])$ conditionally on $\widehat{\gamma}_N$.
This convergence, together with the uniform convergence of $\widehat
{\gamma}_N$ to $\gamma$ presented in Theorem \ref{teo2},
is stronger than the result of Theorem \ref{conditional weak
functional convergence for Gaussian processes} required to build
simultaneous confidence bands.

First, the finite-dimensional convergence of $(\widehat{G}_N)$ to $G$
conditionally on $\widehat{\gamma}_N$
is a trivial consequence of Theorem \ref{teo2}.

Second, we show the tightness of $(\widehat{G}_N)$ in $C([0,T])$
(conditionally on $\widehat{\gamma}_N$)
similarly to the study of $(A_{3,N})$ in the proof of Theorem \ref
{teo2}.
We start by considering the random pseudo-metric
$\hat{d}_{\gamma}^{ m} (s,t) = \mathbb{E}[( \widehat{G}_N(s) -
\widehat{G}_N(t) )^m   |   \widehat{\gamma}_N ] $, where $m\ge1$
is an arbitrary integer.
%Taking into account that $\Delta_{kl} = \Delta_{lk},$ for $k,l =1,
By the moment properties of Gaussian random variables and by (A1), it
holds that
\renewcommand{\theequation}{\arabic{equation}}
\setcounter{equation}{40}
\begin{eqnarray}
\label{random pseudo-metric for th3} \hat{d}_{\gamma}^{ m} (s,t)
& =& C_m \biggl[ \frac{1}{N} \sum_{k,l\in U}
\frac{\Delta_{kl}}{\pi_{kl}} \frac{I_kI_l}{\pi_k \pi_l} \bigl( \widehat{X}_k(s) -
\widehat{X}_k(t) \bigr) \bigl( \widehat {X}_l(s) -
\widehat{X}_l(t) \bigr) \biggr]^{m/2}
\nonumber
\\
& \le& C_m \biggl[ \frac{1}{N} \sum
_{k,l\in U} \frac{|\Delta_{kl}|}{\pi_{kl}} \frac{I_kI_l}{\pi_k \pi_l} \bigl( \widehat
{X}_k(s) - \widehat{X}_k(t) \bigr)^2
\biggr]^{m/2}
\nonumber\\
& \le& C_m \biggl[ \frac{2}{N} \sum
_{k,l\in U} \frac{|\Delta_{kl}|}{\pi_{kl}} \frac{I_kI_l}{\pi_k \pi_l} \bigl(
\tilde{X}_k(s) - \tilde{X}_k(t) \bigr)^2\\
&&\hspace*{17pt}{} +
\frac{2}{N} \sum_{k,l\in U} \frac{|\Delta_{kl}|}{\pi_{kl}}
\frac
{I_kI_l}{\pi_k \pi_l} \bigl( \tilde{\varepsilon}_k(s) - \tilde {
\varepsilon}_k(t) \bigr)^2 \biggr]^{m/2}
\nonumber
\\
& \le& C_m \biggl[ \frac{1}{N} \sum
_{k} \bigl( \tilde{X}_k(s) -
\tilde{X}_k(t) \bigr)^2 \biggr]^{m/2} +
C_m \biggl[ \frac{1}{N}\sum_{k}
\bigl( \tilde{\varepsilon}_k(s) - \tilde{\varepsilon}_k(t)
\bigr)^2 \biggr]^{m/2}.\nonumber
\end{eqnarray}
Note that the the value of the constant $C_m$ varies across the
previous bounds.
Clearly, the first sum in the right-hand side of \eqref{random
pseudo-metric for th3} is dominated by $|s-t|^{m\beta}$ thanks to (A2)
and the approximation properties of local linear smoothers. The second
sum can be viewed as a random quadratic form.
Denoting a square root of $\mathbf{V}_N$ by $\mathbf{V}_N^{1/2}$,
we can write $ \bolds{\varepsilon}_k $ as $\mathbf
{V}_N^{1/2}\mathbf{Z}_k $ for $k=1, \ldots, N$ (the equality
holds in distribution),
where the $\mathbf{Z}_k $ are i.i.d. centered $d$-dimensional
Gaussian vectors with identity covariance matrix.
% matrix $\mathbf{P}_d,$ where $\mathbf{P}_d$ is the projection onto
%the sub-space of $\mathbb{R}^d$ generated by the eigenvectors of $
Thus,
\begin{eqnarray}
\label{random part of pm as qf} \frac{1}{N} \sum
_{k} \bigl( \tilde{\varepsilon}_k(s) - \tilde {
\varepsilon}_k(t) \bigr)^2 & =& \bigl( W(s)-W(t)
\bigr)' \biggl( \frac{1}{N} \sum_{k}
\bolds{\varepsilon}_k \bolds{\varepsilon}_k'
\biggr) \bigl( W(s)-W(t) \bigr)
\nonumber
\\
& \le&\bigl\llVert W(s) - W(t) \bigr\rrVert^2 \biggl\llVert
\frac{1}{N} \sum_{k} \bolds{
\varepsilon}_k \bolds{\varepsilon}_k'
\biggr\rrVert
\\
& \le&\bigl\llVert W(s) - W(t) \bigr\rrVert^2 \llVert
\mathbf{V}_N\rrVert \biggl\llVert \frac{1}{N} \sum
_{k} \mathbf{Z}_k \mathbf{Z}_k'
\biggr\rrVert .\nonumber
\end{eqnarray}
Now, the vector norm $ \llVert  W(s) - W(t) \rrVert^2$ has already
been studied in \eqref{entropy calculation 1} and the sequence
$(\| \mathbf{V}_N \|)$ is bounded by (A4).
The remaining matrix norm in \eqref{random part of pm as qf} is
smaller than the largest eigenvalue, up to a factor $N^{-1},$
of a $d$-variate Wishart matrix with $N$ degrees of freedom.
By (A3) it holds that $d=o(N/ \log\log N)$ and one can apply Theorem
3.1 in Fey \textit{et al.} \cite{FeHoKl08}, which states that for any fixed
$\alpha\ge1,$
\begin{equation}
\label{gdev:fey} \lim_{N \rightarrow\infty} -\frac{1}{N} \log\mathbb{P} \biggl(
\biggl\llVert \frac{1}{N} \sum_{k}
\mathbf{Z}_k \mathbf{Z}_k' \biggr\rrVert
\geq \alpha \biggr)  = \frac{1}{2} ( \alpha- 1 -\log\alpha ).
\end{equation}
An immediate consequence of \eqref{gdev:fey} is that $\llVert  \frac
{1}{N} \sum_{k} \mathbf{Z}_k \mathbf{Z}_k' \rrVert $ remains almost
surely bounded as $N\to\infty$.
Note that the same result holds if instead of (A3), $(d/N)$ remains
bounded away from zero and infinity, thanks to the pioneer work of
Geman \cite{Ge80} on the norm of random matrices.
%and immediately obtain that $ \left\| N^{-1} \sum_{k} \boldsymbol{
%bounded as $N\to\infty$. \hnote{voir aussi apres autre formulation}
Thus, there exists a deterministic constant $C\in(0,\infty)$ such that
\begin{equation}
\hat{d}_{\gamma}^{ m} (s,t) \le C |s-t|^{m\beta} +
\frac
{C}{(dh)^{m/2}} \biggl( \frac{ |s-t|}{h} \wedge1 \biggr)^m
\end{equation}
for all $s,t\in[0,T]$, with probability tending to 1 as $N\to\infty$.
Similarly to the previous entropy calculations,
one can show that there exists a constant $C\in(0,\infty)$ such that
$N(x,\hat{d}_{\gamma}) \le C ( x^{-1/\beta} + (dh^3)^{-1/2} x^{-1} )
$ for all $x \le(dh)^{-1}$
with probability tending to 1 as $N\to\infty$. Applying the maximal
inequality of van der Vaart and Wellner \cite{VaWe00} (Theorem 2.2.4)
to the\vspace*{1pt} conditional increments of $\widehat{G}_N$, with $\phi(t)=t^m$
(usual $L^m$-norm),
one finds a covering integral $\int_0^{1/(ph)}(N(x,\hat{d}_{\gamma
}))^{1/2} \,\mathrm{d}x $ of the order of $(dh)^{ 1/(m\beta) - 1 } +
(dh^3)^{-1/(2m)} (dh)^{1/m-1} $.
Hence, the covering integral tends to zero in probability, provided
that $h\to0 $ and $dh^{(1+1/(2m))/(1-1/(2m))} \to\infty$ as
$N\to\infty$.
Obvisouly, the latter condition on $h$ holds for some integer $m \ge1$
if $dh^{1+\alpha} \to\infty$ for some real $\alpha>0$.
Under this condition, the sequence $(\widehat{G}_N)$ is tight in
$C([0,T])$ %that conditionally on $(\widehat{\gamma}_N)$,
and therefore converges to~$G$.
\end{pf*}
\end{appendix}

%and can be studied with the same arguments used for $A_{3,N}$ in the
%proof of Theorem \ref{teo2}. More
%precisely,
% define the pseudo-metric $d_{\xi}^m((s,t),(s',t')) = \mathbb{E} |
%Then it holds true that
%d_{\xi}&((s,t),(s',t'))\nonumber\\
%& \le C \left\| W(s) - W(t) - W(s' ) + W(t') \right\| \left\| W(s) -
%W(t) + W(s' ) - W(t') \right\| \nonumber\\
%& \le\frac{C}{dh} \left\{ \left( \frac{|s-t|}{h} \wedge1 \right) +

%Hence the diameter of $[0,T]^2$ for $d_{\xi}$ is of order $(dh)^{-1}$.
%Also, for all $0<x<C(dh)^{-1}$,
%the covering number $N(x,d_\xi)$ is no larger than a multiple of
%$(d^2h^4x^2)^{-1}$.
%Therefore, using again the maximal inequality expressed in Corollary
%2.2.5 of van der Vaart and Wellner (2000),
%with the Orlicz norm set as the $L_m$-norm, we obtain that
% \end{equation}
%% and consequently,
%%\begin{equation}
%%\Big\| \sup_{s,t \in[0,T]} | \xi_N(s,t)| \Big\|_m \le\left\|
%% \end{equation}
%%
%%
%% by Markov's inequality,
%%
%%
%This last inequality guarantees the tightness of $(\xi_N)$ in
%$C([0,T]^2)$ provided that $dh^{1+2/m} \to\infty$ as $N\to\infty$.
%This
%in turn holds true for some integer $m\ge1$ if $dh^{1+\alpha} \to
%
%Just as in \eqref{entropy calculation 1}, it can be seen that
%& \le\frac{C}{dh} \left( \frac{|s-t |}{h} \wedge1 \right)^2 .
%Also, by the moment properties of Gaussian random variables, we have
%that
%& \le\frac{C}{N(dh)^2} \left( \frac{|s-t |}{h} \wedge1 \right)^4 .

\section*{Acknowledgements}
 We thank the two anonymous
referees whose careful review lead to substantial improvement of the
paper. Etienne Josserand thanks the Conseil R\'egional de Bourgogne for
its financial support (Faber PhD grant).
\iffalse

%

%suskaldyti doi
\fi

% imsref loaded by imikolaityte, 2012-08-06 14:06:36
%

\printhistory

\end{document}